\documentclass[12pt]{article}
\usepackage[utf8]{inputenc}
\usepackage{amssymb}
\usepackage{amsmath}
\usepackage{amsfonts}
\usepackage{multirow}
\usepackage{tabularx}
\usepackage{makecell}
\usepackage{comment}
\usepackage{float}
\usepackage{graphicx}
\usepackage{listings}
\usepackage{algorithm}
\usepackage{algpseudocode}
\usepackage{pgfplots}
\usepackage{tabularray}
\usepackage{url}
\UseTblrLibrary{booktabs}
\AfterEndEnvironment{lstlisting}{\leavevmode}

\definecolor{metropolis}{HTML}{27373B}
\definecolor{alertcolor}{HTML}{EB811B}
\definecolor{denis}{HTML}{14B03D}
\definecolor{Brown}{cmyk}{0,0.81,1,0.60}
\definecolor{OliveGreen}{cmyk}{0.64,0,0.95,0.40}
\definecolor{CadetBlue}{cmyk}{0.62,0.57,0.23,0}
\definecolor{lightlightgray}{gray}{0.9}
\definecolor{green}{HTML}{44933B}

\usepackage{todonotes}
\presetkeys{todonotes}{inline}{}

\newcommand{\ux}{{\mathbf{x}}}
\newcommand{\uw}{{\mathbf{w}}}
\newcommand{\uc}{{\mathbf{c}}}
\newcommand{\uR}{{\mathbf{R}}}
\newcommand{\uC}{{\mathbf{C}}}
\newcommand{\un}{{\mathbf{n}}}
\newcommand{\ua}{{\mathbf{a}}}

\title{A comparison of Algebraic Multigrid Bidomain solvers on hybrid CPU-GPU architectures}

\author{ Edoardo Centofanti\thanks{
    Dipartimento di Matematica, Università di Pavia,
    Via Adolfo Ferrata, 5, Pavia, 27100, Italy.
    E-mail: {\sf edoardo.centofanti01@universitadipavia.it}.
}  \and Simone Scacchi\thanks{
        Dipartimento di Matematica, Università di Milano,
        Via Cesare Saldini, 50, Milano, 20133, Italy.
        E-mail: {\sf simone.scacchi@unimi.it}.
        }
   }

\date{}
\begin{document}

\maketitle

\begin{abstract}
    The numerical simulation of cardiac electrophysiology is a highly challenging problem in scientific computing.
    The Bidomain system is the most complete mathematical model of cardiac bioelectrical activity.
    It consists of an elliptic and a parabolic partial differential equation (PDE), of reaction-diffusion type,
    describing the spread of electrical excitation in the cardiac tissue.
    The two PDEs are coupled with a stiff system of ordinary differential equations (ODEs), 
    representing ionic currents through the cardiac membrane. 
    Developing efficient and scalable preconditioners for the linear systems arising from the discretization of such computationally challenging model 
    is crucial in order to reduce the computational costs required by the numerical simulations of cardiac electrophysiology. 
    In this work, focusing on the Bidomain system as a model problem, we have benchmarked two popular implementations of the Algebraic Multigrid (AMG) preconditioner embedded in the PETSc library and we have studied the performance on the calibration of specific parameters. 
    We have conducted our analysis on modern HPC architectures, performing scalability tests on multi-core and multi-GPUs setttings. 
    The results have shown that, for our problem, although scalability is verified on CPUs, GPUs are the optimal choice,
    since they yield the best performance in terms of solution time.
    \end{abstract}

\section{Introduction}

Developing physiologically and morphologically realistic and detailed computer models of integrated cardiac function is an important research area,
which may lead to the development of personalised diagnostic techniques and targeted therapies, see e.g. \cite{pot06,aug16,qua17,fed23, jia20, lau23}. 
Also multiscale data available to research can be included in these in-silico models in order to reach such an 
ambitious goal. The main drawback to this approach consists in the associated computational cost, which can rapidly 
become important in terms of solution time of the linear systems deriving from the discretization of the mathematical models
employed to this purpose. 

High-performance computing (HPC) resources are a powerful research tool in order to partially fulfill the problem, 
but there is still a huge amount of effort to optimize the resources employed and 
make a sustainable use of them. In particular, over the last decade, General-Purpose Graphic Processing Units (GPGPUs) have shown an extraordinary computational power
with respect to using a large number of CPUs to solve large numerical systems. This has led to GPUs being a significant item in the budgets of 
companies and institutions interested in assembling an HPC cluster. 

Apart from computational resources, remarkable efforts have been made also in the development of efficient numerical schemes and techniques
aimed at reducing the number of iterations of an iterative method employed to solve a large linear system. This is the case of the preconditioners, which can be tailored according to the problem to be addressed or to be part of a larger class of general purpose preconditioning techniques. 
The Algebraic Multigrid (AMG) \cite{run87} is both a solver and a general purpose preconditioner. Due to its versatility, it can be employed as an effective preconditioner for reducing Conjugate Gradient (CG) iterations when solving particularly expensive systems arising from the discretization of elliptic partial differential equations (PDEs), which are largely employed in the mathematical modeling of cardiac electrophysiology.

In this work, we focused our efforts on solving efficiently the cardiac Bidomain model, a macroscopic representation of the cardiac tissue modeling the spatio-temporal evolution of the intra- and extracellular electric potentials. 
In the formulation considered throughout the work, it consists of two PDEs, an elliptic and a parabolic one, coupled with a system of ordinary differential equations (ODEs) representing the ionic currents through the cellular membrane. The model homogenizes both the intra- and the extracellular space, namely in each physical degree of freedom (DOF) both spaces cohexist. Previous work on solvers involving AMG for the Bidomain system are \cite{pen09, pen11, pla07}, but other preconditioners have been successfully tested and employed, such as geometric multigrid and domain decomposition preconditioners, as Multilevel Schwarz\cite{pav08, pav11}, Neumann Neumann\cite{zam13,zam14} and BDDC\cite{huy22-1, huy22-2}.

The novelty of this work consists in having chosen and tested two different implementations of the AMG preconditioner and having benchmarked their performances on an HPC architecture involving both CPUs and GPUs, in order to verify if GPUs are effective for this particular problem and setup. A previous work on this topic \cite{nei12} has shown a speedup of roughly a factor 10 when employing GPUs with respect to the best performance on CPU, but the benchmark is limited up to 20 GPUs and the Hypre implementation of AMG considered, BoomerAMG, was not yet implemented in order to run on GPUs. In this work we will also expand some of those results exploiting our \textit{in-house} codebase, which considers a slightly different numerical scheme.

The rest of the paper is organized as follows:
in Section \ref{model} the Bidomain model is presented and the numerical schemes employed are explained and commented; in Sections \ref{parset} and \ref{arch}
we present the parameters, the simulation setup and the HPC architectures considered; in Section \ref{sec:amg} we describe the AMG implementations and the importance
of two different threshold parameters employed; in Sections \ref{teststruct} and \ref{testunstr} we present the numerical tests performed and the corresponding results.

\section{Model}\label{model}

The model studied in the following is the macroscopic Bidomain model of electrocardiology in the parabolic-elliptic formulation\cite{fra14,ven09,pen05}. 
Denoting by $\Omega\subset\mathbb{R}^3$ the bounded physical region occupied by the portion of myocardium of our interest, we represent the tissue as the superimposition
of two anisotropic continuous media, called intra- and extracellular media. In this model, it is assumed that they coexist at every point
and are separated by a distributed continuous cellular membrane.  
Our problem consists of finding the extracellular potential 
$u_{e}: \Omega\times (0,T)\rightarrow \mathbb{R}$,
the transmembrane potential $v: \Omega\times (0,T)\rightarrow \mathbb{R}$ and the gating and ionic concentration variables
$\uw: \Omega\times (0,T)\rightarrow \mathbb{R}^{N_w}$  and $\uc: \Omega\times (0,T)\rightarrow \mathbb{R}^{N_c}$, respectively, such that: 
\begin{equation}
\left\{    
\begin{array}{ll}
\displaystyle \chi C_m \frac{\partial v}{\partial v} - \text{div}(D_i\nabla v) - \text{div}(D_i\nabla u_e) + I_{ion} (v,\uw,\uc) = I_{app}^i &  \text{in } \Omega\times (0,T),\vspace{0.2cm}\\
\displaystyle -\text{div}(D_i\nabla v) -\text{div} ((D_i + D_e)\nabla u_e) = I_{app}^i + I_{app}^e & \text{in } \Omega\times (0,T),\vspace{0.2cm}\\
\displaystyle \frac{\partial \uw}{\partial t} - \uR(v,\uw) = 0 & \text{in } \Omega\times (0,T),\vspace{0.2cm}\\
\displaystyle \frac{\partial \uc}{\partial t} - \uC(v,\uw,\uc) = 0 & \text{in } \Omega\times (0,T),\vspace{0.2cm}\\
\displaystyle \un^\top D_i\nabla (v + u_e) = 0 & \text{in } \partial\Omega\times (0,T),\vspace{0.2cm}\\
\displaystyle \un^\top D_e\nabla u_e = 0 & \text{in } \partial\Omega\times (0,T),\vspace{0.2cm}\\
\displaystyle v(\ux,0) = v_0(\ux), \quad \uw(\ux,0)=\uw_0(\ux), \quad \uc(\ux,0)=\uc_0(\ux) & \text{in } \Omega, 
\end{array}        
\right.
\end{equation}
where $I_{ion}$ is the ionic current related to the ionic model. In this work we have employed the ten Tusscher-Panfilov\cite{ten06} model (TP06). The functions $\uR$ and $\uC$ describe the dynamics of the gating and ionic concentration variables $\uw$ and $\uc$, respectively. $I_{app}^i$ and $I_{app}^e$ denote the applied currents in the intra- and extracellular region, respectively, that must satisfy the compatibility condition
\[
\int_\Omega (I_{app}^i+I_{app}^e)\,d\ux=0.
\] 
The anisotropy of the tissue, due to the fiber arrangement of cardiac myocytes (see e.g. \cite{leg95}), is described by the conductivity tensors $D_{i,e}(\ux)$ at any point $\ux\in\Omega$, which have the following definition:
\begin{equation}
    \begin{split}
        D_{i,e}(\ux) & = \sigma_{l}^{i,e} \ua_l(\ux)\ua_l^\top(\ux) + \sigma_{t}^{i,e} \ua_t(\ux)\ua_t^\top(\ux) + \sigma_{n}^{i,e} \ua_n(\ux)\ua_n^\top(\ux) = \\
                            & = \sigma_{t}^{i,e} I + (\sigma_l^{i,e}-\sigma_t^{i,e}) \ua_l(\ux)\ua_l^\top(\ux) + (\sigma_n^{i,e}-\sigma_t^{i,e}) \ua_n(\ux)\ua_n^\top(\ux), 
    \end{split}
\end{equation}
where $\ua_l$, $\ua_t$ and $\ua_n$ represent a triple of orthonormal principal axes, with $\ua_l$ parallel 
to the local fiber direction, $\ua_n$ and $\ua_t$ orthogonal and tangent to the radial laminae, respectively, and both transversal to the fiber axis. 
Furthermore, we have denoted as $\sigma_l^{i,e}$, $\sigma_t^{i,e}$ and $\sigma_n^{i,e}$ the conductivity coefficients for 
the intra- and extracellular media along the corresponding directions $\ua_l$, $\ua_t$ and $\ua_n$. 
We recall that these tensors model the structure of the cardiac tissue, i.e.
an ensemble of fibers rotating counterclockwise from epi- (the outermost protective layer of the heart) to endocardium (the innermost), organized as
muscle foils running radially from epi- to endocardium. 

\subsection{Variational formulation}

Let us define $V:=H^1(\Omega)$ the usual Sobolev space and then define 
\[\widetilde{V}=\{\psi\in V: \int_\Omega \psi = 0\},\] \[U=V\times\widetilde{V}=\{u=(\phi,\psi):\phi\in V, \psi \in \widetilde{V}\},\]
the $L^2$-inner product $(\phi,\psi)=\int_{\Omega}\phi\psi d\ux$ $\forall \phi,\psi\in L^2(\Omega)$, and the elliptic bilinear forms
\[a_{i,e}(\phi,\psi)=\int_\Omega (\nabla\phi)^\top D_{i,e}(\ux)\nabla\psi d\ux,\]
\[a(\phi,\psi) = \int_\Omega (\nabla\phi)^\top D(\ux)\nabla\psi d\ux,\quad \forall\phi,\psi\in H^1(\Omega),\]
where $D=D_i+D_e$ is the bulk conductivity tensor.

The variational formulation of the Bidomain model reads as follows:
given $v_0\in L^2(\Omega),\,\uw_0\in L^2(\Omega)^{N_w},\,\uc_0\in L^2(\Omega)^{N_c}$, $I^{i,e}_{app}\in L^2(\Omega\times (0,T))$, find $v\in L^2(0,T;V)$, $u_e\in L^2(0,T;\widetilde{V})$, $\uw\in L^2(0,T;L^2(\Omega)^{N_w})$ and $\uc\in L^2(0,T;L^2(\Omega)^{N_c})$ such that $\frac{\partial v}{\partial t}\in L^2(0,T;V)$,
$\frac{\partial \uw}{\partial t}\in L^2(0,T;L^2(\Omega)^{N_w})$, $\frac{\partial \uc}{\partial t}\in L^2(0,T;L^2(\Omega)^{N_c})$ and $\forall t\in (0,T)$
\[
\left\{
\begin{array}{ll}
\displaystyle c_m\frac{\partial}{\partial t}(v,\hat{v})+a_i(v+u_e,\hat{v})+(I_{ion}(v,\uw,\uc),\hat{v})=(I^i_{app},\hat{v}) & \forall\hat{v}\in V\vspace{0.2cm} \\
\displaystyle a_i(v,\hat{u}_e)+a(u_e,\hat{u}_e)=0 & \forall \hat{u}_e\in\widehat{V} \vspace{0.2cm}\\
\displaystyle\frac{\partial}{\partial t}(w_j,\hat{w})-(R_j(v,\uw),\hat{w})=0, & \forall\hat{w}\in V,\,j=1,...,N_w \vspace{0.2cm}\\
\displaystyle\frac{\partial}{\partial t}(c_j,\hat{c})-(C_j(v,\uw,\uc),\hat{c})=0, & \forall\hat{c}\in V,\,j=1,...,N_c
\end{array}
\right.
\]
with initial conditions $v=v_0$, $\uw=\uw_0$ and $\uc=\uc_0$ and where we have imposed $I^e_{app}=-I^i_{app}$.

\subsection{Space and Time Discretization}

Let $\mathcal{T}_h$ be a quasi-uniform triangulation of $\Omega$ having maximal diameter $h$
and $V_h$ be an associate conforming finite element space. Let us select a finite element basis $\{\phi_p\}_{p=1}^N$, evaluate $M=\{m_{pj}\}$, the diagonal mass matrix, with the usual mass-lumping technique and 
$A_{i,e}=\{a_{pj}^{i,e}\}$ the symmetric intra- and extracellular stiffness matrices, having elements
\[a_{pj}^{i,e} = \int_{\Omega} D_{i,e}\nabla\phi_j\cdot\nabla\phi_p dx\]
By means of a standard Galerkin procedure, we can rewrite the semi-discrete Bidomain problem, discretized in space, in the following compact form:
\begin{equation}\label{eq:bido-matrix}
    \begin{cases}
        \displaystyle c_m\mathcal{M}\frac{d}{dt}\begin{bmatrix}
            \mathbf{v} \\ \mathbf{u}_e
        \end{bmatrix} + \mathcal{A} \begin{bmatrix}
            \mathbf{v} \\ \mathbf{u}_e \end{bmatrix} + \begin{bmatrix}
                M I_{ion}(\mathbf{v},\mathbf{w},\mathbf{c}) \\ \mathbf{0}
            \end{bmatrix} = \begin{bmatrix}
                M\mathbf{I}_{app}^i \\ \mathbf{0}
            \end{bmatrix} \\
            \displaystyle \frac{d\mathbf{w}}{dt} = \mathbf{R}(\mathbf{v},\mathbf{w}) \vspace{0.2cm}\\
            \displaystyle \frac{d\mathbf{c}}{dt} = \mathbf{C}(\mathbf{v},\mathbf{w},\mathbf{c})
    \end{cases}    
\end{equation}
where we denote the block mass and stiffness matrices as
\[\mathcal{M}=\begin{bmatrix} 
M & 0 \\ 0 & 0 
\end{bmatrix}%
\qquad\hfill%
\mathcal{A}=\begin{bmatrix}
    A_i & A_i \\ A_i & A_i+A_e
\end{bmatrix}\]
and we define $\mathbf{v}$, $\mathbf{u}_e$, $\mathbf{w}=(\mathbf{w}_1,...,\mathbf{w}_{N_w})^\top$,
$\mathbf{c}=(\mathbf{c}_1,...,\mathbf{c}_{N_c})^\top$,
$R(\mathbf{v},\mathbf{w})=(R_1(\mathbf{v},\mathbf{w}),...,R_{N_w}(\mathbf{v},\mathbf{w}))^\top$,
$C(\mathbf{v},\mathbf{w},\mathbf{c})=(C_1(\mathbf{v},\mathbf{w},\mathbf{c}),...,C_{N_c}(\mathbf{v},\mathbf{w},\mathbf{c}))^\top$,
$\mathbf{I}_{ion}(\mathbf{v},\mathbf{w})$ and $\mathbf{I}_{app}^e$  as the coefficient vectors of 
finite element approximations of $u_i$, $u_e$, $v$, $w_r$, $R_r(v,w_1,...,w_{N_w})$, $C_r(v,w_1,...,w_{N_w},c_1,...,c_{N_c})$, $I_{ion}(v,w_1,...,w_{N_w})$
and $I_{app}^{i,e}$ respectively. \\

For the time discretization, we consider an implicit-explicit (IMEX) strategy, which consists of decoupling 
the ODEs from the PDEs and of treating the linear diffusion terms implicitly and the non-linear reaction terms 
explicitly. 
We then solve uncoupled the two equations (\ref{eq:bido-matrix}) in discretized form. \\
In particular, considering $\mathbf{v}^n$ and $\mathbf{u}_e^n$ at the timestep $n$, 
we solve the elliptic equation evaluating $\mathbf{u}_e^{n+1}$, and then we solve the parabolic equation in order to update 
the transmembrane potential $\mathbf{v}^{n+1}$. We summarize the scheme, given $\mathbf{w}^n$,   
$\mathbf{v}^n$, $\mathbf{u}_e^n$, as

\begin{align}
    &\mathbf{w}^{n+1} + \Delta t \mathbf{R}(\mathbf{v}^n,\mathbf{w}^{n+1}) &&= \mathbf{w}^n \\
    & \mathbf{c}^{n+1} + \Delta t \mathbf{C}(\mathbf{v}^n,\mathbf{w}^{n+1},\mathbf{c}^{n+1}) &&= \mathbf{c}^n \\
    &(A_i+A_e)\mathbf{u}_e^n &&= -A_i \mathbf{v}^n  \\
    &\left(\frac{c_m}{\Delta t}M+A_i\right)\mathbf{v}^{n+1} &&= \frac{c_m}{\Delta t}M\mathbf{v}^n-A_i\mathbf{u}_e^n + M\mathbf{I}_{ion}(\mathbf{v}^n,\mathbf{w}^{n+1},\mathbf{c}^{n+1}) +M\mathbf{I}_{app}^{i,n}.
\end{align}

Overall, at each timestep we solve once the linear system with matrix $A_i+A_e$, i.e. the discrete form of the elliptic equation,
and once the linear system with matrix $\frac{c_m}{\Delta t}M+A_i$ deriving from the parabolic equation. Being the resulting 
systems very large, due to the number of DOFs considered, they must be solved with an iterative method. We have chosen in particular the Preconditioned Conjugate Gradient (PCG) method, 
since the matrix arising from the parabolic equation is symmetric positive definite whereas that arising from the elliptic
equation is symmetric positive semi-definite. The preconditioners used for the parabolic and elliptic systems will be discussed in the next sections.

\section{Parameters and Setting}\label{parset}
We have focused our study on two geometries: a truncated ellipsoid, modeling an idealized left ventricle (LV) and a realistic
geometry, representing a patient specific LV at three different resolutions. The idealized LV is discretized by a structured hexaedral mesh, whereas the patient specific LV is discretized by an unstructured mesh consisting of irregular hexaedra.

The truncated ellipsoid is built following the parametric equations:

\begin{equation}
    \begin{cases}
        & x = a(r)\cos\theta\cos\phi,\quad \hfill \theta_{min}\le\theta\le\theta_{max} \\
        & y = b(r)\cos\theta\sin\phi, \hfill \phi_{min}\le\phi\le\phi_{max} \\
        & z = c(r)\sin\phi, \hfill 0\le r\le 1, 
    \end{cases}    
\end{equation}
where $a(r)= a_1 + r(a_2-a_1)$, $b(r)=b_1 + r(b2-b1)$ and $c(r)= c_1 + r(c2-c1)$, with $a_{1,2}$, $b_{1,2}$, $c_{1,2}$
coefficients defining the main axes of the ellipsoid. The geometric parameters are reported in Table \ref{geomp}.

For what concerns the patient specific LV geometry, denoted in the following as ``U-mesh" we have considered three different refinements, which are reported in table \ref{zygtab}. The mesh has been provided us by Marco Fedele at Mox Laboratory, Politecnico di Milano, and the fibers have been
generated using the open-source code lifex-fiber \cite{afr23}, developed at Mox Laboratory, Politecnico di Milano; see also \cite{pie21,fed23}. 

The ionic membrane model considered is the ten Tusscher-Panfilov model (TP06)\cite{ten06}. The stimulus is applied for $1 \text{ ms}$ with intensity of $350 \text{ mA}/\text{cm}^3$. Depending on the numerical tests performed, the simulation spans the first 5ms of the excitation phase as well as 500 ms, corresponding to almost a full cycle, both with timestep of 0.05 ms.

\begin{table}
    \centering
\begin{tblr}{cc|cc}
    \toprule
    Parameter & Value & Parameter & Value \\
    \midrule
    $a_1$ & 2.2 & $\theta_{min}$ & $-\frac{3\pi}{2}$ \\ 
    $a_2$ & 3.3 & $\theta_{max}$ & $\frac{\pi}{2}$ \\
    $b_1$ & 2.2 & $\phi_{min}$ & $-\frac{3\pi}{8}$ \\
    $b_2$ & 3.3 & $\phi_{max}$ & $\frac{\pi}{8}$ \\
    $c_1$ & 5.9 & & \\ 
    $c_2$ & 6.4 & & \\
    \bottomrule
\end{tblr}
\caption{Geometric parameters of idealized LV domain}
\label{geomp}
\end{table}

\begin{table}
\centering
\begin{tblr}{c}
  \toprule
  Name & Physical DOFs & Elements \\
  \midrule
  U-mesh 1 & 35725 & 30108 \\ 
  U-mesh 2 & 258415 & 240864 \\
  U-mesh 3 & 1987285 & 1926912 \\
  \bottomrule
\end{tblr}
\caption{Number of points and elements for the ``U-mesh" geometry}
\label{zygtab}
\end{table}

\section{Algebraic Multigrid}\label{sec:amg}
We are interested in solving a problem of the form
\begin{equation}\label{sys}
    A\mathbf{x} = \mathbf{f},
\end{equation}
with $A\in\mathbb{R}^{n\times n}$, $\mathbf{x}$, $\mathbf{f}$ $\in \mathbb{R}^n$. \\
Being this kind of systems usually very large, they cannot be solved with direct methods in short time. \\
Thus it is necessary to employ iterative methods and other strategies in order to get a fast,
albeit approximate solution. \\
The main idea of the Algebraic Multigrid (AMG) algorithm is solving (\ref{sys}) by cycling through levels 
composed of coarse (i.e. smaller) linear systems and finding updates that, projected on the original space,
improve the solution. In this way, the so called \lq smooth error\rq $\,e$ 
that is not eliminated through iterative relaxations, is removed
by coarse-grid correction. \\
The idea is implemented by solving the residual equation $Ae=r$
on a coarser grid, then through interpolation the solution is brought back to the 
finer grid and finally the fine-grid approximation is updated $u \longleftarrow u+e$. \\
Given the matrix $A$ with entries $a_{ij}$ and the various indices $i,j$ identified as grid points 
on a grid $\Omega$, the actors playing a role in the AMG are the following:
\begin{enumerate}
    \item $M$ grids, i.e. index sets $\Omega =\Omega^1\supset \Omega^2 \supset ... \supset \Omega^M$.
    \item $M$ grid operators $A^1,A^2,...,A^M$.
    \item $M-1$ Interpolation operators $P^1,P^2,...,P^{M-1}$.
    \item $M-1$ Restriction operators $R^1,R^2,...,R^{M-1}$.
    \item $M-1$ Smoothers $S^1,S^2,...,S^{M-1}$.
\end{enumerate}
All these components are defined in the \textit{setup phase} of the algorithm:

\begin{algorithm}
\caption{AMG setup phase}
\begin{algorithmic}
    \For{$k=1;\,\, k< M;\,\, k++$}
    \Statex $\quad$ \text{Partition $\Omega^k$ into disjoint sets $C^k$ and $F^k$.}
    \Statex $\quad$ \text{Set $\Omega^{k+1} = C^k$.}
    \Statex $\quad$ \text{Define interpolation $P^k$.}
    \Statex $\quad$ \text{Define restriction $R^k$ (often $R^k=(P^k)^\top$).}
    \State $A^{k+1} \gets R^kA^kP^k$
    \Statex $\quad$ \text{Set up $S^k$.}
    \EndFor
\end{algorithmic}
\end{algorithm}
When the setup phase is completed, the algorithm proceeds with a recursively defined cycle.
\newpage
Following the notation in \cite{ste06}, we call this phase `MGV',
since it is often addressed as `Multigrid V-Cycle'. Other cycles (W and F) are possible,
but in this work we will focus only on V. The steps are the following:

\begin{algorithm}
\caption{MGV($A^k,R^k,P^k,S^k,u^k,f^k$)}
\begin{algorithmic}
    \If{$k==M$}
    \Statex \text{$\quad$ solve $A^Mu^M=f^M$ with a direct solver.}
    \Else
    \Statex \text{$\quad$ apply the smoother $S^k$ $\mu_1$ times to $A^ku^k=f^k$.}
    \State $r^k\gets f^k-A^ku^k$ \Comment{Coarse grid correction step}
    \State $r^{k+1}\gets R^kr^k$
    \Statex \text{$\quad$ apply MGV($A^{k+1},R^{k+1},P^{k+1},S^{k+1},e^{k+1},r^{k+1}$)} \Comment{Recursive step}
    \State $e^k \gets P^ke^{k+1}$ \Comment{Interpolation step}
    \State $u^k\gets u^k+e^k$ \Comment{Correction}
    \Statex \text{$\quad$ apply the smoother $S^k$ $\mu_2$ times to $A^ku^k=f^k$.}
    \EndIf
\end{algorithmic}
\end{algorithm}

We note that the smoothing step is performed through a Richardson iteration of the form $u^k_{j+1} = u^k_j + \omega (S^k)^{-1} (f-Au^k_j)$,
with $S^k$ usually Jacobi, Gauss-Seidel or ILU and $\omega$ relaxation factor.\\
Throughout this work, we will employ two popular AMG implementations in the field of HPC and many-core systems: GAMG, which is the built-in AMG solver in the PETSc library\cite{pet23} and BoomerAMG, provided within the Hypre library\cite{fal06}. They are both contained in the PETSc library, the Portable, Extensible Toolkit for Scientific Computation, which includes a large suite of scalable parallel linear and nonlinear equation solvers, ODE integrators, and optimization algorithms\cite{pet23}.

A critical step consists of the construction of the restriction matrices $R^k$, which are involved in the coarsening phase. \\
Since the algorithm should not rely directly on the geometry of the domain, many algebraic techniques have been explored and developed over the years.
Although many coarsening algorithms are available for the implementations employed \cite{yan02,gri06a,gri06b,stu99}, we have considered the ones suggested as default options which are, for GAMG, PETSc default AMG implementation, a modified
Maximal Independent Set (MIS) algorithm \cite{ada98,ada04}, while for Hypre a Hybrid Maximal Independent Set (HMIS) algorithm \cite{des06}. \\
In the first phase of the modified MIS, a graph is built from the nodes $i,j$ of the elements $a_{i,j}$ of the matrix $A=A^1$ and a weigth $w_{ij} = \frac{a_{ij}}{\sqrt{a_{ii}a_{jj}}}$ is attributed to each edge $(i,j)$. 
A threshold on the edge weigth is thus set such that at each coarsening step all the edges with weigth less that the threshold are cut. Then, a greedy MIS or the parallel implementation (Luby) is applied to the modified graph. We call $S$ the resulting set. 
Clusters $C_j$ are defined with the following procedure: for each $i\notin S$, $j\in S$, $i\in C_j$ if $w_{ij}=\max_{\bar{j}\in S}w_{i\bar{j}}$. Finally, prolongators for AMG are defined as 
\[ P_{ij}=\begin{cases} |C_j|^{-\frac{1}{2}} \text{ if }i\in C_j \\ 0 \text{ otherwise}\end{cases} \]
In the first phase of HMIS instead, $A$ is explored and for each row the coarsening nodes are chosen between the ones satisfying the condition
\begin{equation}
    |a_{i,j}|\ge\alpha \max_{k\neq i}|a_{i,k}|
\end{equation}  
with $\alpha$ called \textit{strong threshold} parameter. This brings to the definition of sets of \textit{strongly connected nodes} which are used to define the sets $C_j$\\
In the following sections we will explore how the threshold and the strong threshold influence our code performances.

\section{Architectures and Technical Setup}\label{arch}
Except where otherwise stated, the numerical experiments have been performed on the Marconi100
cluster at CINECA laboratory. 
Marconi100 is a Linux Cluster with 980 nodes, each one including 2 x 16 cores 
IBM POWER9 AC922 @ 3.1 GHz and 4 NVIDIA Volta V100 GPUs with 16GB memory per GPU 
and NVlink 2.0. Each node has a memory of 256GB. 
Computations have been performed using up to 16 nodes considering either CPU or GPU
performance. The maximum number of cores was 512, while we used up to 64 GPUs. 
Our \textit{in-house} code is written in C with CUDA kernels for solving the membrane model on the GPUs. 
The numerical aspects of the code are based on the PETSc library, while the communications between parallel 
processors is handled by MPI.
Tests have been focused mainly on the performance and the behaviour at different scales
of the preconditioners used for the elliptic equation of the model.
Therefore, we have tested two implementations of the AMG preconditioner, available in PETSc:
\begin{itemize}
    \item The GAMG preconditioner, the default implementation in PETSc, on CPU (up to 512, 16 nodes). We have also
    performed a few tests on GPU (up to 4, 1 node), but our current setup did not allow to perform an exaustive 
    analysis on GPUs. In particular we have observed out of memory issues and a general suboptimal
    use of the device (GPU) resources with our implementation.
    \item The Hypre BoomerAMG preconditioner, on CPU (up to 512, 16 nodes) and on GPU (up to 64, 16 nodes).
\end{itemize}
Regarding the parabolic equation in the formulation considered, we have employed Block Jacobi preconditioner on CPU, while we have left the system unpreconditioned on GPU.
For the PETSc options employed in this work, refer to appendix \ref{opts}.
We have divided the results into subsections: For the structured mesh, in \textit{Test 1} we tune the AMG 
threshold hyperparameter both for Hypre BoomerAMG and GAMG, in \textit{Test 2} we have performed a strong
scaling test, fixing a global dimension for the problem while changing the size of the CPUs and GPUs employed, while
in \textit{Test 3} we have performed a weak scaling test,
keeping the local size (i.e. the size per processor or per GPU) of the problem unchanged, while varying the number of 
processors or GPUs. \\
For the unstructured mesh, in \textit{Test 1} we have tuned the threshold like in the structured case, 
while in \textit{Test 2} we have performed a strong scaling test, fixing the global size of the problem 
while varying the number of processors or GPUs for the three refinements of the mesh. \\

\section{Numerical results on structured meshes}\label{teststruct}
We have first considered the structured mesh, discretizing the truncated ellipsoid described in section \ref{model}.
\subsection{Test 1 - AMG Threshold calibration}
In the first test, we have studied the behaviour of the average solution time per timestep for the elliptic 
equation while varying the threshold parameters, i.e. the threshold for the GAMG preconditioner and the strong 
threshold for Hypre. These parameters act in the coarsening step of the algorithm and the best choices for them are problem dependent\cite{pet23,fal06}.
For each implementation of the preconditioner employed, we have performed tests with different threshold parameters on a single node, exploiting 4 CPUs or 4 GPUs for each test. The geometry considered is the truncated ellipsoid discretized by a $32\times 32\times 16$ hexaedral mesh.
For each multigrid implementation we have reported the CG iterations and the solution times for the elliptic system.
In Tables \ref{threshhg}, \ref{threshhc} and \ref{threshgc} we have reported the results for Hypre on GPU, Hypre on CPU and for GAMG on CPU, respectively.
As expected, only CG iterations and solution times for the elliptic equation are affected for varying the threshold parameters. Optimal results are obtained with a strong threshold of 0.4 for Hypre on GPU, 0.5 for Hypre on CPU and 0.06 for GAMG.
\begin{table}[H]
    \centering
    \begin{tblr}{columns={c,font=\footnotesize},vline{2}}
    \toprule
    Threshold & 0.25 & 0.3 & 0.4 & 0.5 & 0.6 & 0.7 \\
    $\text{It}_{\text{ellip, mean}}$ & 23.30 & 19.29 & 13.92 & 22.81 & 24.57 & 32.11 \\
    $T_{\text{ellip, mean}}$ (s) & 9.6E-02 & 7.4E-02 & 5.6E-02 & 0.10 & 0.11 & 0.15 \\
    \bottomrule
\end{tblr}
\caption{AMG Threshold calibration, structured mesh. Results for Hypre GPU solver. $\text{It}_{\text{ellip, mean}}$: average CG iterations per timestep for the elliptic solver. $T_{\text{ellip, mean}}$: average CG solution time (in s) per timestep for the elliptic solver.}
\label{threshhg}
\end{table}

\begin{table}[H]
    \centering
    \begin{tblr}{columns={c,font=\footnotesize},vline{2}}
    \toprule
    Threshold & 0.25 & 0.3 & 0.4 & 0.5 & 0.6 & 0.7 \\
    \midrule
    $\text{It}_{\text{ellip, mean}}$ & 6.84 & 7.76 & 11.90 & 6.27 & 8.84 & 11.30 \\
    $T_{\text{ellip, mean}}$ (s) & 2.9E-02 & 3.1E-02 & 4.5E-02 & 2.4E-02 & 3.1E-02 & 3.6E-02 \\
    \bottomrule
    \end{tblr}
\caption{AMG Threshold calibration, structured mesh. Results for Hypre CPU solver. Same format as in Table \ref{threshhg}.}
\label{threshhc}
\end{table}

\begin{table}[H]
\centering
\begin{tblr}{columns={c,font=\footnotesize},vline{2}}
    \toprule
    Threshold & 0.0 & 0.01 & 0.02 & 0.03 & 0.04 & 0.05 & 0.06 & 0.07 \\
    \midrule
    $\text{It}_{\text{ellip, mean}}$ & 143.1 & 66.17 & 64.03 & 52.47 & 44.85 & 42.18 & 9.50 & 9.64 \\
    $T_{\text{ellip, mean}}$ (s) & 0.44 & 0.16 & 0.17 & 0.18 & 0.18 & 0.18 & 0.06 & 0.07 \\
    \bottomrule
\end{tblr}
\caption{AMG Threshold calibration, structured mesh. Results for GAMG CPU solver. Same format as in Table \ref{threshhg}.}
\label{threshgc}
\end{table}

\subsection{Test 2 - strong scaling}
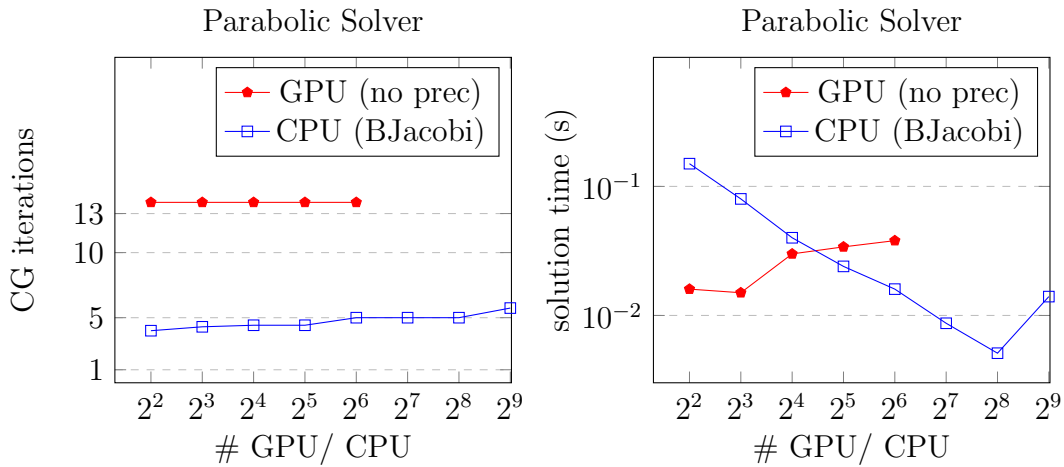
\begin{figure}[H] 
    \begin{tikzpicture}
        \def\titleax{Parabolic Solver}
        \def\titlex{\# GPU/ CPU}
        \def\titley{CG iterations}
        \def\xmin{0}
        \def\xmax{520}
        \def\ymin{0}
        \def\ymax{25}
        \def\legpos{north east}
        \def\wannagridsx{false}
        \def\wannagridsy{true}
        
        \def\coordinatesxya{(4,13.86)(8,13.86)(16,13.86)(32,13.86)(64,13.86)}
        \def\legenda{GPU (no prec)} 
       
        \def\coordinatesxyb{(4,4.00)(8,4.30)(16,4.42)(32,4.42)(64,5.00)(128,5.00)(256,5.00)(512,5.74)}
        \def\legendb{CPU (BJacobi)}
    
        \begin{axis}[
            width=0.5\textwidth,
            title={\titleax},
            xlabel={\titlex},
            ylabel={\titley},
            xmin=\xmin, xmax=\xmax,
            ymin=\ymin, ymax=\ymax,
            xtick={4,8,16,32,64,128,256,512},
            ytick={1,5,10,13},
            legend pos=\legpos,
            xmajorgrids=\wannagridsx,
            ymajorgrids=\wannagridsy,
            grid style=dashed,
            xmode=log,
            log basis x=2
        ]
    
        \addplot[
            color=red,
            mark=pentagon*,
            ]
            coordinates {
            \coordinatesxya
            };
    
        \addplot[
            color=blue,
            mark=square,
            ]
            coordinates {
            \coordinatesxyb
            };
            
        \legend{\legenda,\legendb}
    
        \end{axis}
    \end{tikzpicture}%
    \begin{tikzpicture}
        \def\titleax{Parabolic Solver}
        \def\titlex{\# GPU/ CPU}
        \def\titley{solution time (s)}
        \def\xmin{0}
        \def\xmax{520}
        \def\ymin{0}
        \def\ymax{1}
        \def\legpos{north east}
        \def\wannagridsx{false}
        \def\wannagridsy{true}
        
        \def\coordinatesxya{(4,1.6e-02)(8,1.5e-02)(16,3.0e-02)(32,3.4e-02)(64,3.8e-02)}
        \def\legenda{GPU (no prec)} 
    
        \def\coordinatesxyb{(4,0.15)(8,8e-02)(16,4.0e-02)(32,2.4e-2)(64,1.6e-2)(128,8.7e-3)(256,5.1e-3)(512,1.4e-2)}
        \def\legendb{CPU (BJacobi)}
    
        \begin{axis}[
            width=0.5\textwidth,
            title={\titleax},
            xlabel={\titlex},
            ylabel={\titley},
            xmin=\xmin, xmax=\xmax,
            ymin=\ymin, ymax=\ymax,
            xtick={4,8,16,32,64,128,256,512},
            ytick={1e-4,1e-3,1e-2,1e-1},
            legend pos=\legpos,
            xmajorgrids=\wannagridsx,
            ymajorgrids=\wannagridsy,
            grid style=dashed,
            xmode=log,
            log basis x=2,
            ymode=log,
            log basis y=10
        ]
    
        \addplot[
            color=red,
            mark=pentagon*,
            ]
            coordinates {
            \coordinatesxya
            };

        \addplot[
            color=blue,
            mark=square,
            ]
            coordinates {
            \coordinatesxyb
            };
            
        \legend{\legenda,\legendb}
    
        \end{axis}
        \end{tikzpicture}
        \caption{Strong scaling, structured meshes. Comparison of CG iterations (left) and solution time (right) for the parabolic solver on GPU or CPU.}
        \label{strong1}
    \end{figure}

\begin{figure}[H] 
    \begin{tikzpicture}
        \def\titleax{Membrane Model Solver}
        \def\titlex{\# GPU/ CPU}
        \def\titley{solution time (s)}
        \def\xmin{0}
        \def\xmax{520}
        \def\ymin{0}
        \def\ymax{1.1}
        \def\legpos{north east}
        \def\wannagridsx{false}
        \def\wannagridsy{true}
        
        \def\coordinatesxya{(4,8.4043e-04)(8,5.1412e-04)(16,3.3960e-04)(32,1.9872e-04)(64,1.8497e-04)}
        \def\legenda{GPU} 
       
        \def\coordinatesxyb{(4,0.5089)(8,0.2587)(16,0.1294)(32,0.0707)(64,0.0367)(128,0.0183)(256,0.0102)(512,0.0051)}
        \def\legendb{CPU}
    
        \begin{axis}[
            width=0.5\textwidth,
            title={\titleax},
            xlabel={\titlex},
            ylabel={\titley},
            xmin=\xmin, xmax=\xmax,
            ymin=\ymin, ymax=\ymax,
            xtick={4,8,16,32,64,128,256,512},
            ytick={1e-4,1e-3,1e-2,1e-1},
            legend pos=\legpos,
            xmajorgrids=\wannagridsx,
            ymajorgrids=\wannagridsy,
            grid style=dashed,
            xmode=log,
            log basis x=2,
            ymode=log,
            log basis y=10
        ]
    
        \addplot[
            color=green,
            mark=halfcircle*,
            ]
            coordinates {
            \coordinatesxya
            };
    
        \addplot[
            color=red,
            mark=pentagon*,
            ]
            coordinates {
            \coordinatesxyb
            };
        \legend{\legenda,\legendb}
    
        \end{axis}
        \end{tikzpicture}%
        \begin{tikzpicture}
        \def\titleax{Total solution time per timestep}
        \def\titlex{\# GPU/ CPU}
        \def\titley{solution time (s)}
        \def\xmin{0}
        \def\xmax{520}
        \def\ymin{0}
        \def\ymax{1300}
        \def\legpos{north east}
        \def\wannagridsx{false}
        \def\wannagridsy{true}
        
        \def\coordinatesxya{(4,7.1e-2)(8,9.7e-2)(16,0.4)(32,0.22)(64,0.29)}
        \def\legenda{GPU hypre} 
       
        \def\coordinatesxyb{(4,3.60)(8,1.06)(16,0.74)(32,0.36)(64,0.48)(128,0.23)(256,0.35)(512,1.05)}
        \def\legendb{CPU hypre}
    
        \def\coordinatesxyc{(4,18.17)(8,11.13)(16,5.77)(32,3.60)(64,2.34)(128,1.33)(256,0.92)(512,1.01)}
        \def\legendc{CPU GAMG}
    
        \begin{axis}[
            width=0.5\textwidth,
            title={\titleax},
            xlabel={\titlex},
            ylabel={\titley},
            xmin=\xmin, xmax=\xmax,
            ymin=\ymin, ymax=\ymax,
            xtick={4,8,16,32,64,128,256,512},
            ytick={1e-4,1e-3,1e-2,1e-1,1,10},
            legend pos=\legpos,
            xmajorgrids=\wannagridsx,
            ymajorgrids=\wannagridsy,
            grid style=dashed,
            xmode=log,
            log basis x=2,
            ymode=log,
            log basis y=10
        ]
    
        \addplot[
            color=green,
            mark=halfcircle*,
            ]
            coordinates {
            \coordinatesxya
            };
    
        \addplot[
            color=red,
            mark=pentagon*,
            ]
            coordinates {
            \coordinatesxyb
            };
    
        \addplot[
            color=blue,
            mark=square,
            ]
            coordinates {
            \coordinatesxyc
            };
            
        \legend{\legenda,\legendb,\legendc}
    
        \end{axis}
        \end{tikzpicture}
        \caption{Strong scaling, structured meshes. Comparison of solution time for the membrane model on GPU or CPU (left). Total solution time per timestep, given as the sum of parabolic, elliptic and membrane solution times on GPU or CPU (right).}
        \label{strong2}
    \end{figure}
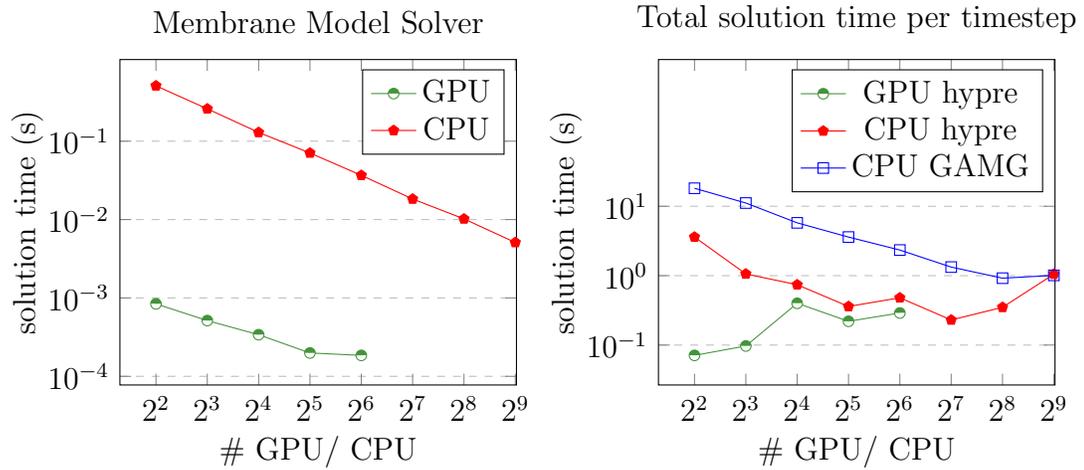

We have performed here a strong scaling test, by fixing the global size of the problem while increasing the number of GPUs or CPUs. As in the previous test, we have studied the average solution time per timestep for the elliptic and parabolic equations and for the membrane model (TP06). We have also reported the average number of CG iterations per time step for both the parabolic and elliptic systems. Regarding the elliptic system, the setup of AMG is the one described in the previous section. We fix the global size of the problem by employing a $128\times 128\times 64$ mesh, leading
to a total of 2163330 DOFs.
The results reported in figures \ref{strong1} and \ref{strong2} have shown that all solvers are scalable in terms of CG iterations, which remain almost bounded when increasing the number of the CPUs or GPUs. The solution times on GPU are not scalable, whereas on CPU they are scalable up to 128/256 cores. However, the solution time for the elliptic system with the Hypre GPU solver is significantly lower than with the CPU Hypre and GAMG solvers.
The results reported in Figure \ref{strong2} (left panel) show that the membrane model solver is scalable both on CPU and GPU. Moreover, the solution times on GPU are about two order of magnitude lower than on CPU. 
Figure \ref{strong2} (right panel) shows that, in terms of total solution time, the best performance is obtained for GPU. This performance cannot be achieved using CPU.

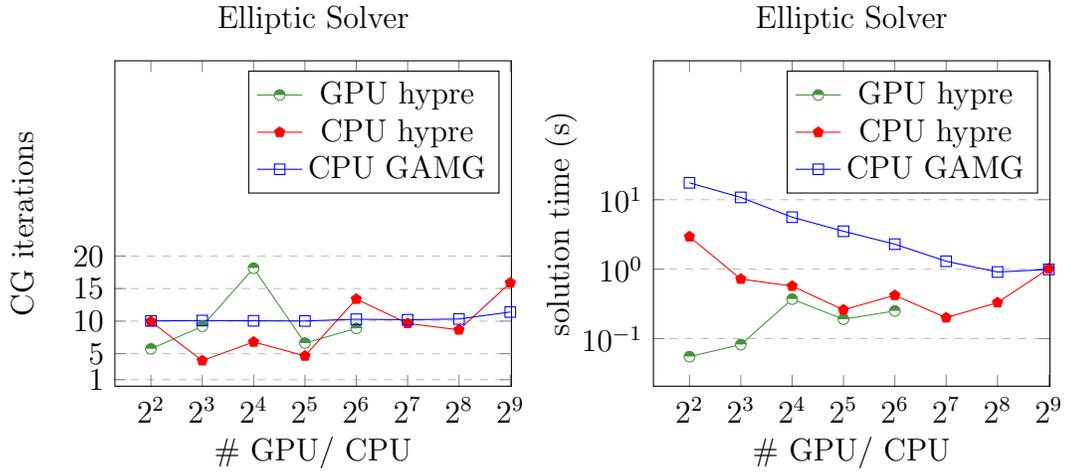
\begin{figure}[H] 
    \begin{tikzpicture}
        \def\titleax{Elliptic Solver}
        \def\titlex{\# GPU/ CPU}
        \def\titley{CG iterations}
        \def\xmin{0}
        \def\xmax{520}
        \def\ymin{0}
        \def\ymax{50}
        \def\legpos{north east}
        \def\wannagridsx{false}
        \def\wannagridsy{true}
        
        \def\coordinatesxya{(4,5.74)(8,9.20)(16,18.14)(32,6.63)(64,8.88)}
        \def\legenda{GPU hypre} 
       
        \def\coordinatesxyb{(4,9.94)(8,3.93)(16,6.80)(32,4.66)(64,13.39)(128,9.65)(256,8.66)(512,15.89)}
        \def\legendb{CPU hypre}
    
        \def\coordinatesxyc{(4,10.03)(8,10.09)(16,10.05)(32,10.02)(64,10.29)(128,10.21)(256,10.33)(512,11.40)}
        \def\legendc{CPU GAMG}
    
        \begin{axis}[
            width=0.5\textwidth,
            title={\titleax},
            xlabel={\titlex},
            ylabel={\titley},
            xmin=\xmin, xmax=\xmax,
            ymin=\ymin, ymax=\ymax,
            xtick={4,8,16,32,64,128,256,512},
            ytick={1,5,10,15,20},
            legend pos=\legpos,
            xmajorgrids=\wannagridsx,
            ymajorgrids=\wannagridsy,
            grid style=dashed,
            xmode=log,
            log basis x=2
        ]
    
        \addplot[
            color=green,
            mark=halfcircle*,
            ]
            coordinates {
            \coordinatesxya
            };
    
        \addplot[
            color=red,
            mark=pentagon*,
            ]
            coordinates {
            \coordinatesxyb
            };
    
        \addplot[
            color=blue,
            mark=square,
            ]
            coordinates {
            \coordinatesxyc
            };
            
        \legend{\legenda,\legendb,\legendc}
    
        \end{axis}
        \end{tikzpicture}%
    \begin{tikzpicture}
        \def\titleax{Elliptic Solver}
        \def\titlex{\# GPU/ CPU}
        \def\titley{solution time (s)}
        \def\xmin{0}
        \def\xmax{520}
        \def\ymin{0}
        \def\ymax{1000}
        \def\legpos{north east}
        \def\wannagridsx{false}
        \def\wannagridsy{true}
        
        \def\coordinatesxya{(4,5.5E-02)(8,8.2E-02)(16,0.37)(32,0.19)(64,0.25)}
        \def\legenda{GPU hypre} 
       
        \def\coordinatesxyb{(4,2.94)(8,0.72)(16,0.57)(32,0.26)(64,0.42)(128,0.20)(256,0.33)(512,1.02)}
        \def\legendb{CPU hypre}
    
        \def\coordinatesxyc{(4,17.51)(8,10.79)(16,5.60)(32,3.51)(64,2.29)(128,1.30)(256,0.91)(512,0.99)}
        \def\legendc{CPU GAMG}
    
        \begin{axis}[
            width=0.5\textwidth,
            title={\titleax},
            xlabel={\titlex},
            ylabel={\titley},
            xmin=\xmin, xmax=\xmax,
            ymin=\ymin, ymax=\ymax,
            xtick={4,8,16,32,64,128,256,512},
            ytick={1e-4,1e-3,1e-2,1e-1,1,10},
            legend pos=\legpos,
            xmajorgrids=\wannagridsx,
            ymajorgrids=\wannagridsy,
            grid style=dashed,
            xmode=log,
            log basis x=2,
            ymode=log,
            log basis y=10
        ]
    
        \addplot[
            color=green,
            mark=halfcircle*,
            ]
            coordinates {
            \coordinatesxya
            };
    
        \addplot[
            color=red,
            mark=pentagon*,
            ]
            coordinates {
            \coordinatesxyb
            };
    
        \addplot[
            color=blue,
            mark=square,
            ]
            coordinates {
            \coordinatesxyc
            };
            
        \legend{\legenda,\legendb,\legendc}
    
        \end{axis}
        \end{tikzpicture}
        \caption{Strong scaling, structured meshes. Comparison of CG iterations (left) and solution time (right) for the elliptic solver on GPU or CPU.}
    \end{figure}

\subsection{Test 3 - Weak scaling}
We have performed here a weak scaling test, by fixing the local size of the problem while increasing the number of GPUs or CPUs. In Tables \ref{weakhg2}-\ref{weakgc2} we have reported the average solution time per timestep for the elliptic and parabolic equations and for the membrane model (TP06). We have studied also the average number of CG iterations per time step for both the parabolic and elliptic systems. Regarding the elliptic system, the setup of AMG is the same previously described. The local size of the problem is fixed such that each CPU/ GPU handles a local mesh of $16 \times 16 \times 16$ elements, for a total amount of 9826 DOFs per worker. 

Data reported in Tables \ref{weakhg2}, \ref{weakhc} and \ref{weakgc2} have shown that in terms of CG iterations the solvers are scalable, with an argument similar to the previous section. Solution times for the elliptic equation with Hypre on GPU are not scalable, but the results seem to be affected by a synchronization overhead due to idle threads in the GPUs. Data relative to CPU, both for GAMG and Hypre show instead good weak scalability, with solution times which are overall comparable, increasing the number of cores, while keeping fixed the local size of the problem. Similar comments can be made regarding the parabolic equation for both CG iterations and solution time. We have also observed very good scalability with both CPUs and GPUs for the solution time of the membrane model. 
\begin{table}[H]
\begin{tblr}{columns={c,font=\footnotesize}}
    \toprule
    Num GPU & DOFs &$\text{It}_{\text{parab, mean}}$& $\text{It}_{\text{ellip, mean}}$ & $T_{\text{memb, mean}}$ (s) & $T_{\text{parab, mean}}$ (s) & $T_{\text{ellip, mean}}$ (s) \\
    \midrule
    4 & 37026 & 25.80 & 11.10 & 1.7E-04 & 1.4E-02 & 5.0E-02 \\
    8 & 71874 & 23.04 & 16.32 & 2.3E-04 & 5.1E-02 & 0.34 \\
    16 & 141570 & 21.09 & 15.82 & 5.2E-04 & 9.3E-02 & 0.70 \\
    32 & 278850 & 17.21 & 7.29 & 1.0E-03 & 0.15 & 0.68 \\
    64 & 549250 & 16.63 & 19.63 & 1.1E-03 & 0.15 & 1.81 \\
    \bottomrule
\end{tblr}
\caption{Weak Scaling on GPU, structured meshes. CG preconditioned by Hypre BoomerAMG GPU as elliptic solver and unpreconditioned CG as parabolic solver. Membrane model solved on GPU. BoomerAMG Threshold = 0.4. $\text{It}_{\text{parab, mean}}$: average CG iterations per timestep for the parabolic solver. $\text{It}_{\text{ellip, mean}}$: average CG iterations per timestep for the elliptic solver. 
$T_{\text{memb, mean}}$: average solution time (in s) per timestep for the membrane model.
$T_{\text{parab, mean}}$: average CG solution time (in s) per timestep for the parabolic solver. $T_{\text{ellip, mean}}$: average CG solution time (in s) per timestep for the elliptic solver.}
\label{weakhg2}
\end{table}

\begin{table}[H]
\begin{tblr}{columns={c,font=\footnotesize}}
    \toprule
    Num CPU & DOFs & $\text{It}_{\text{parab, mean}}$ & $\text{It}_{\text{ellip, mean}}$ & $T_{\text{memb, mean}}$ (s) & $T_{\text{parab, mean}}$ (s) & $T_{\text{ellip, mean}}$ (s) \\
    \midrule
    4 & 37026 & 3.00 & 7.45 & 1.7E-03 & 2.0E-03 & 2.7E-02 \\
    8 & 71874 & 3.00 & 69.55 & 3.1E-03 & 2.1E-03 & 0.31 \\
    16 & 141570 & 3.00 & 36.93 & 6.7E-03 & 2.3E-03 & 0.32 \\
    32 & 278850 & 3.00 & 12.02 & 1.1E-02 & 2.8E-03 & 9.1E-02 \\
    64 & 549250 & 4.90 & 13.84 & 1.4E-02 & 4.5E-03 & 0.13 \\
    128 & 1090050 & 4.97 & 20.48 & 1.4E-02 & 5.2E-03 & 0.25 \\
    256 & 2163330 & 5.00 & 10.59 & 1.4E-02 & 6.7E-03 & 0.22 \\
    512 & 4293378 & 7.87 & 7.77 & 1.3E-02 & 1.1E-02 & 0.21 \\
    \bottomrule
\end{tblr}
\caption{Weak Scaling on CPU, structured meshes. CG preconditioned by Hypre BoomerAMG CPU as elliptic solver and CG preconditioned by Block Jacobi as parabolic solver. Membrane model solved on GPU. BoomerAMG Threshold 0.5. Same format as in Table \ref{weakhg2}.}
\label{weakhc}
\end{table}

\begin{table}[H]
\begin{tblr}{columns={c,font=\footnotesize}}
    \toprule
    Num CPU & DOFs & $\text{It}_{\text{parab, mean}}$ & $\text{It}_{\text{ellip, mean}}$ & $T_{\text{memb, mean}}$ (s) & $T_{\text{parab, mean}}$ (s) & $T_{\text{ellip, mean}}$ (s) \\
    \midrule
    4 & 37026 & 3.00 & 9.50 & 1.9E-04 & 2.1E-03 & 5.9E-02 \\
    8 & 71874 & 3.00 & 10.16 & 7.8E-04 & 2.2E-03 & 0.18 \\
    16 & 141570 & 3.00 & 10.16 & 1.5E-03 & 2.3E-03 & 0.16 \\
    32 & 278850 & 3.00 & 10.27 & 2.9E-03 & 2.9E-03 & 0.16 \\
    64 & 549250 & 4.90 & 10.59 & 3.1E-03 & 4.2E-03 & 0.78 \\
    128 & 1090050 & 4.97 & 10.36 & 3.5E-03 & 7.8E-03 & 0.77 \\
    256 & 2163330 & 5.00 & 10.52 & 3.4E-03 & 1.1E-02 & 0.67 \\
    512 & 4293378 & 7.87 & 13.81 & 3.2E-03 & 2.0E-02 & 4.02 \\
    \bottomrule
\end{tblr}
\caption{Weak Scaling on CPU, structured meshes. CG preconditioned by PETSc GAMG CPU as elliptic solver and CG preconditioned by Block Jacobi as parabolic solver. Membrane model solved on GPU. GAMG Threshold 0.5. Same format as in Table \ref{weakhg2}.}
\label{weakgc2}
\end{table}

\begin{figure}[H]
    \centering
    \begin{tblr}{colspec = {X[c] X[c]}, colsep = 2pt}
    \includegraphics[width=0.5\linewidth]{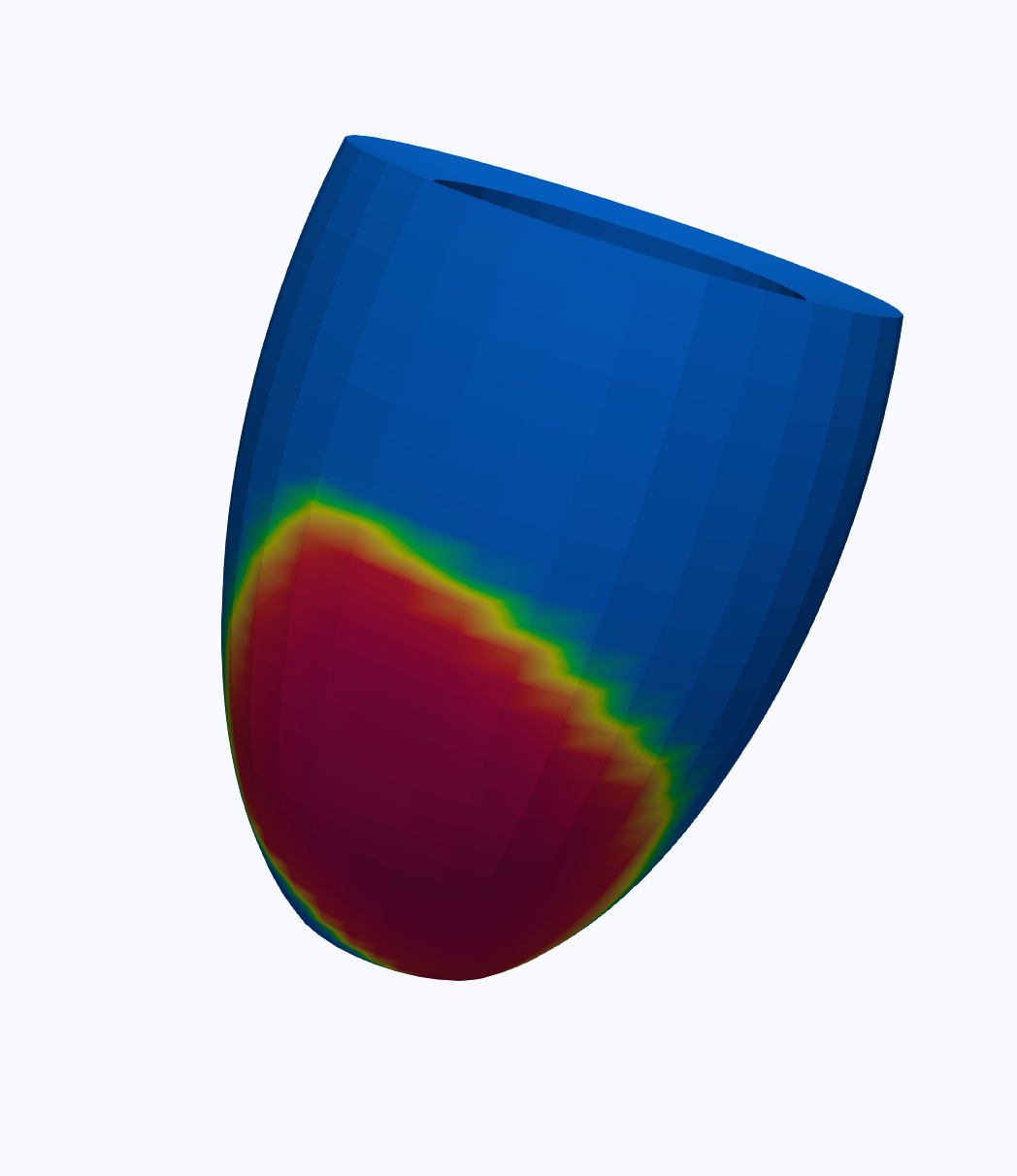} & \includegraphics[width=0.5\linewidth]{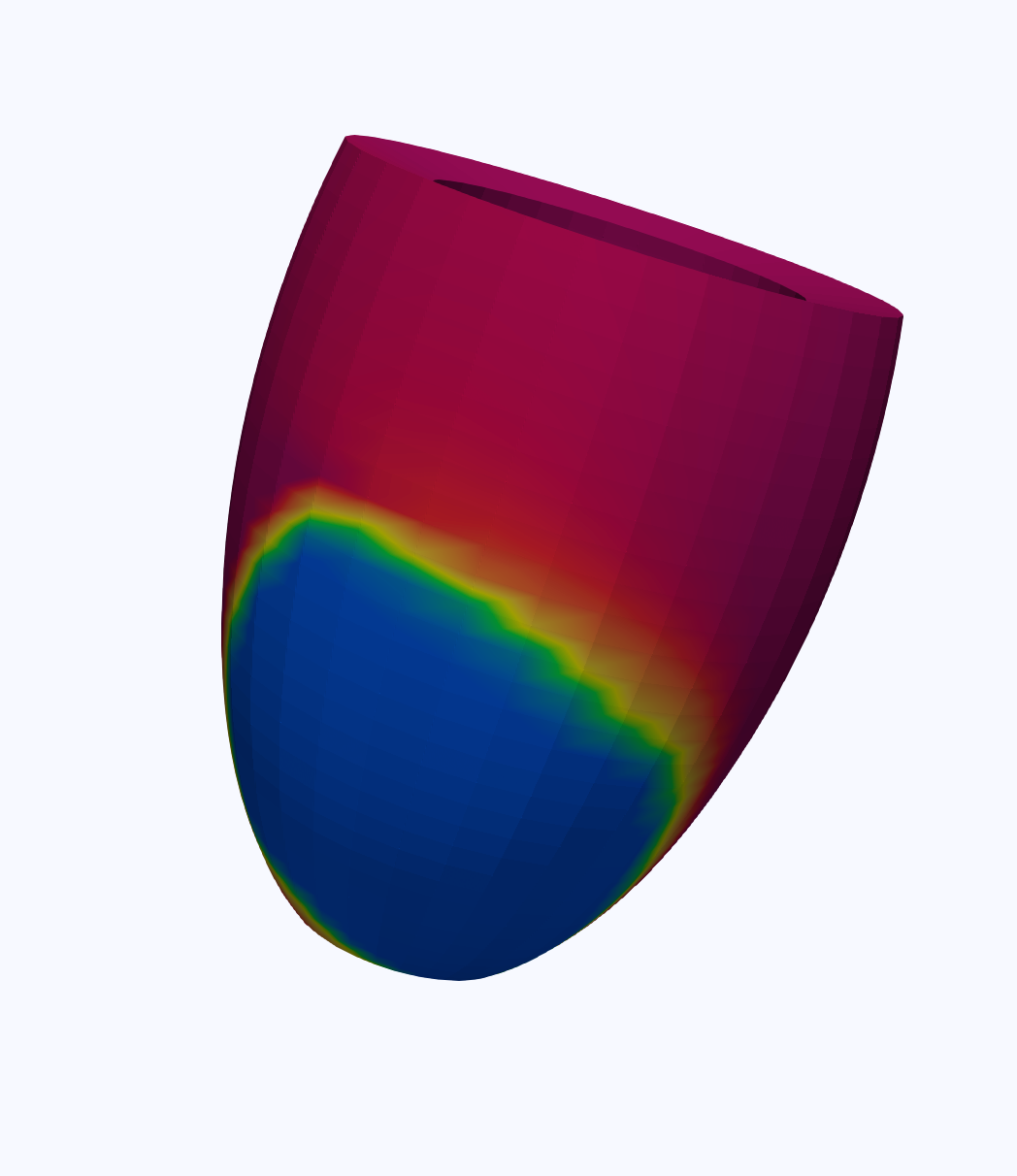} \\
    \includegraphics[width=0.5\linewidth]{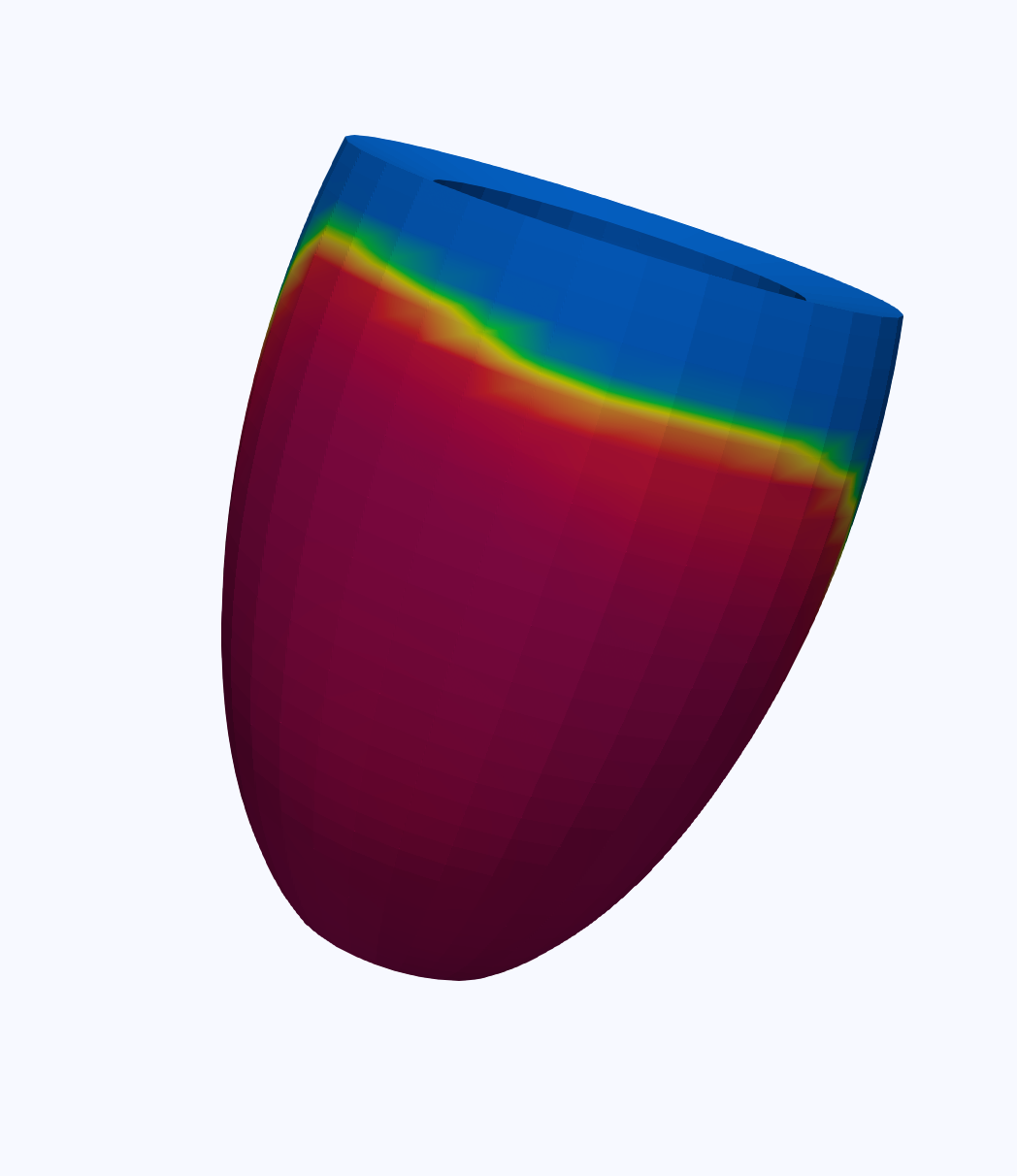} & \includegraphics[width=0.5\linewidth]{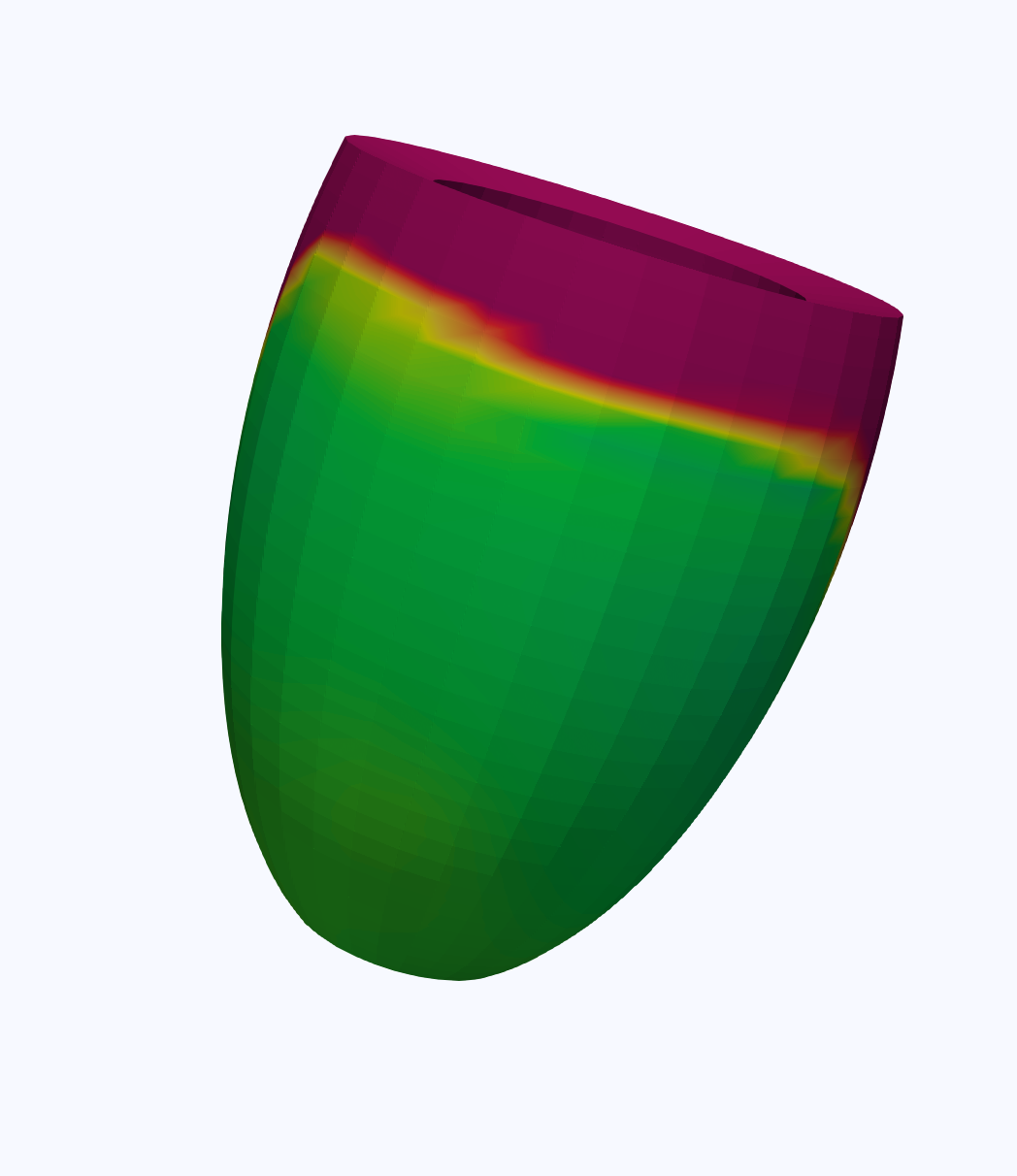} \\
    \includegraphics[width=0.5\linewidth]{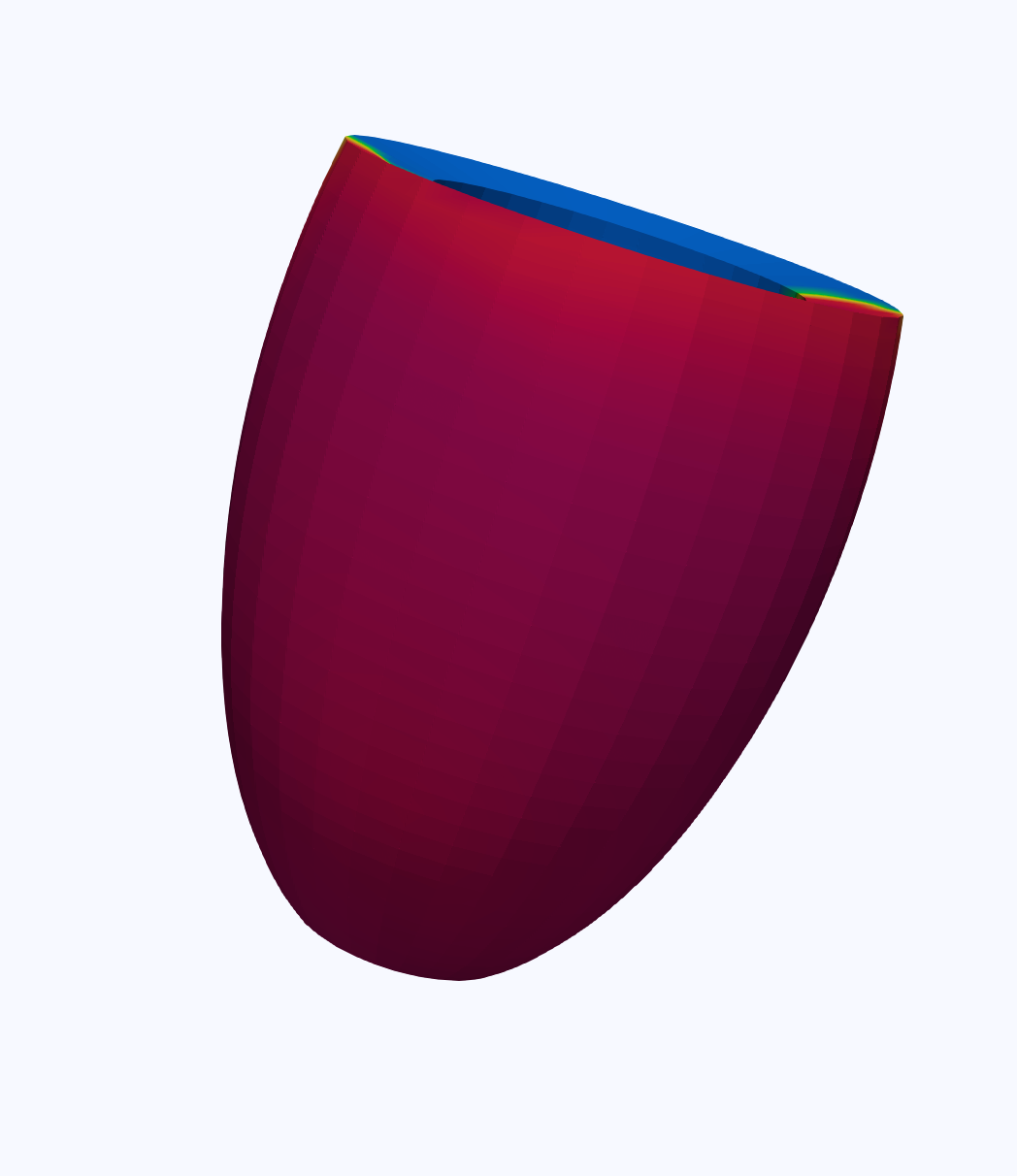} & \includegraphics[width=0.5\linewidth]{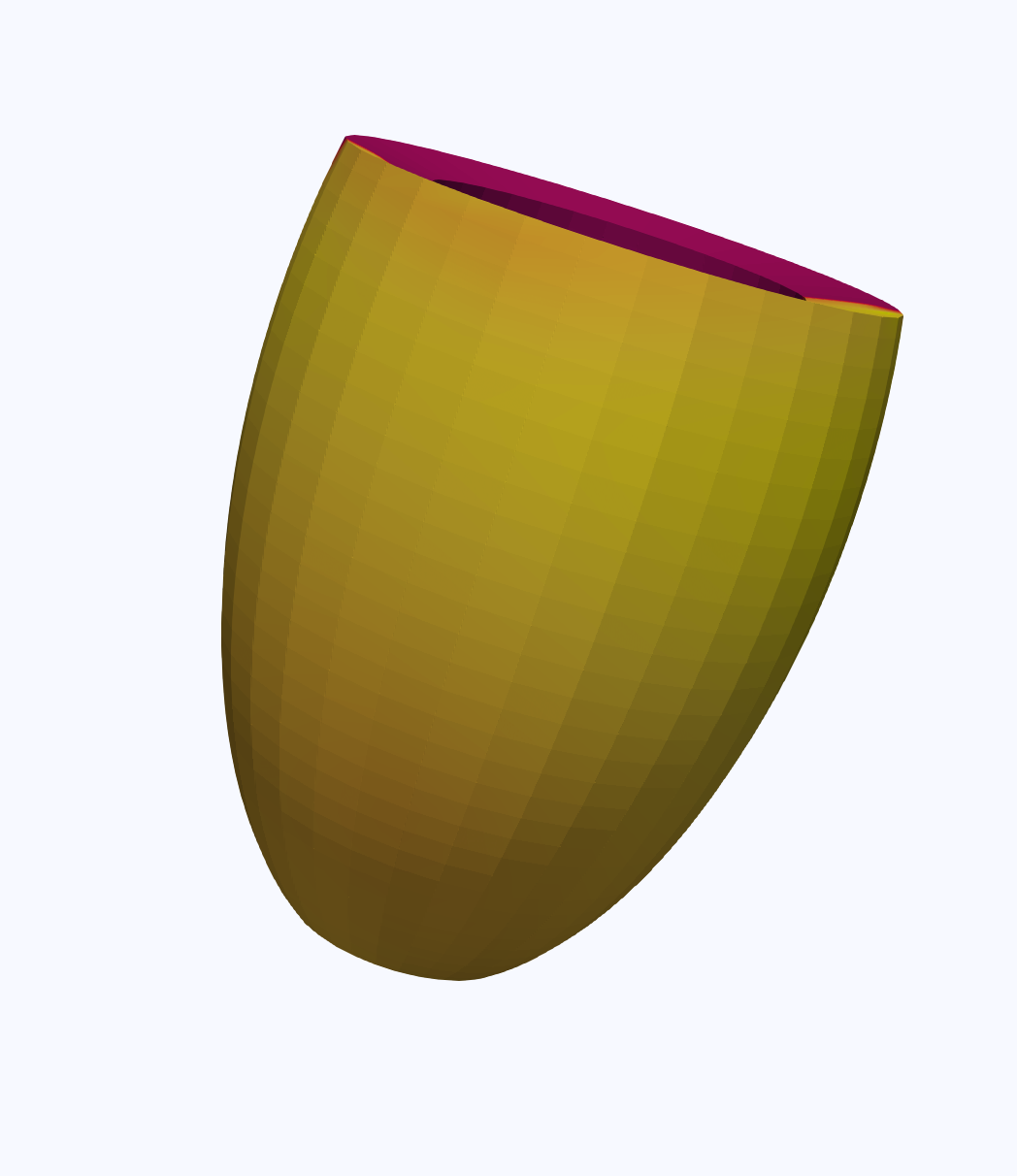} \\
    \includegraphics[width=0.5\linewidth]{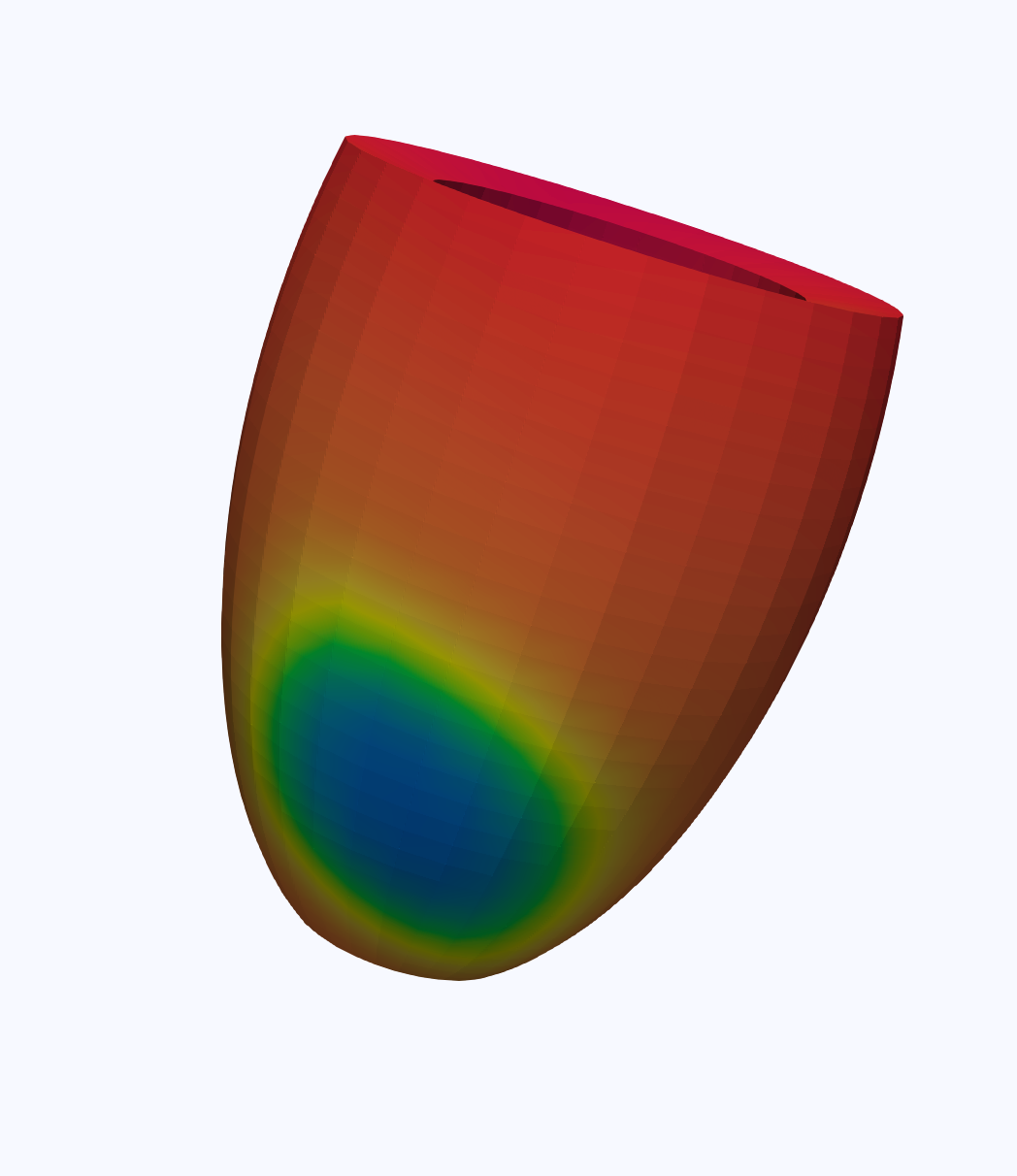} & \includegraphics[width=0.5\linewidth]{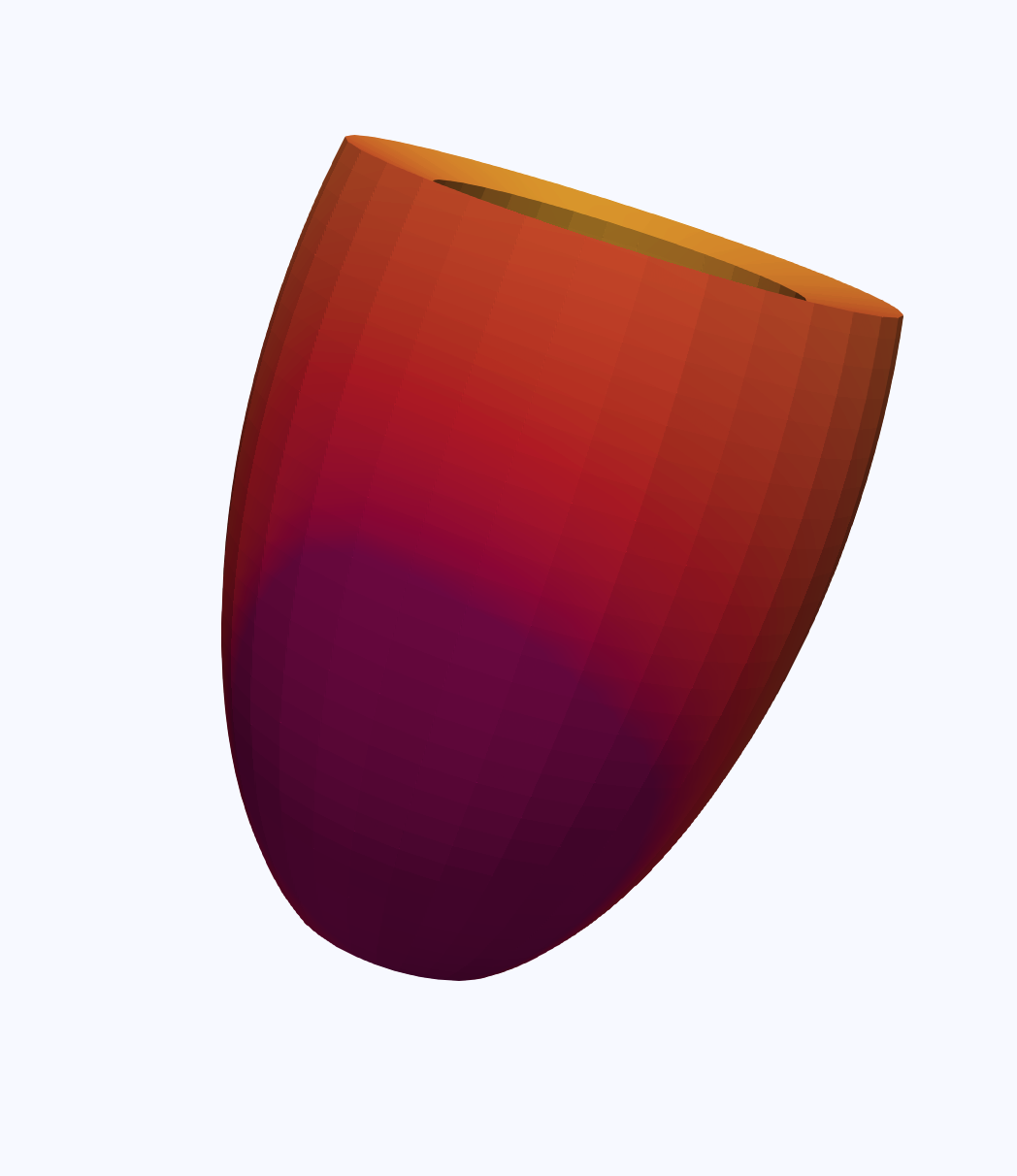} \\
    \includegraphics[width=0.5\linewidth]{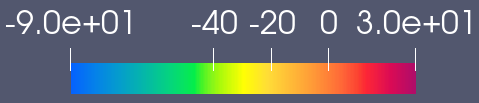} & \includegraphics[width=0.5\linewidth]{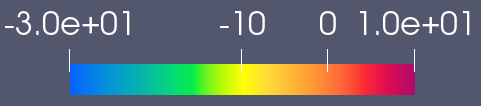} \\
    \end{tblr}
\caption{Transmembrane potential (left) and extracellular potential (right) snapshots on the epicardial surface, represented by the structured mesh. The values of the displayed map are in mV}
\end{figure}

\begin{figure}[H] 
    \begin{tikzpicture}
        \def\titleax{Elliptic Solvers}
        \def\titlex{\# CPU}
        \def\titley{CG iterations}
        \def\xmin{0}
        \def\xmax{520}
        \def\ymin{0}
        \def\ymax{75}
        \def\legpos{north east}
        \def\wannagridsx{false}
        \def\wannagridsy{true}
             
        \def\coordinatesxya{(4,9.4950)(8,10.1584)(16,10.1584)(32,10.2673)(64,10.5941)(128,10.3564)(256,10.5248)(512,13.8119)}
        \def\legenda{gamg}
          
        \def\coordinatesxyb{(4,7.4455e+00)(8,6.9545e+01)(16,3.6931e+01)(32,1.2020e+01)(64,1.3842e+01)(128,2.0475e+01)(256,1.0594e+01)(512,7.7723e+00)}
        \def\legendb{hypre}
    
        \begin{axis}[
            width=0.5\textwidth,
            title={\titleax},
            xlabel={\titlex},
            ylabel={\titley},
            xmin=\xmin, xmax=\xmax,
            ymin=\ymin, ymax=\ymax,
            xtick={4,8,16,32,64,128,256,512},
            ytick={10,20,40,60},
            legend pos=\legpos,
            legend style={font=\tiny},
            xmajorgrids=\wannagridsx,
            ymajorgrids=\wannagridsy,
            grid style=dashed,
            xmode=log,
            log basis x=2
        ]
    
        \addplot[
            color=blue,
            mark=square,
            ]
            coordinates {
            \coordinatesxya
            };
    
            \addplot[
                color=red,
                mark=pentagon*,
                ]
                coordinates {
                \coordinatesxyb
                };
            \legend{\legenda,\legendb}
    
        \end{axis}
    \end{tikzpicture}
    \begin{tikzpicture}
        \def\titleax{Elliptic and Parabolic Solvers}
        \def\titlex{\# CPU}
        \def\titley{solution time (s)}
        \def\xmin{0}
        \def\xmax{520}
        \def\ymin{0}
        \def\ymax{300}
        \def\legpos{north west}
        \def\wannagridsx{false}
        \def\wannagridsy{true}
        \def\xticks{4,8,16,32,64,128,256,512}
        \def\yticks{1e-3,1e-2,1e-1,1,10}
        
        \def\coordinatesxya{(4,0.0591)(8,0.1767)(16,1.1640)(32,2.1598)(64,5.7817)(128,0.7672)(256,0.6669)(512,4.0168)}
        \def\legenda{gamg}
        
        \def\coordinatesxyb{(4,2.7497e-02)(8,3.0859e-01)(16,3.2400e-01)(32,9.1421e-02)(64,1.2883e-01)(128,2.5445e-01)(256,2.2019e-01)(512,2.1097e-01)}
        \def\legendb{hypre}
               
        \def\coordinatesxyc{(4,2.0270e-03)(8,2.1206e-03)(16,2.2937e-03)(32,2.7983e-03)(64,4.5372e-03)(128,5.2048e-03)(256,6.7476e-03)(512,1.0922e-02)}
        \def\legendc{parabolic}
    
        \begin{axis}[
            width=0.5\textwidth,
            title={\titleax},
            xlabel={\titlex},
            ylabel={\titley},
            xmin=\xmin, xmax=\xmax,
            ymin=\ymin, ymax=\ymax,
            xtick={4,8,16,32,64,128,256,512},
            ytick={1e-3,1e-2,1e-1,1,10},
            legend pos=\legpos,
            legend style={font=\tiny},
            xmajorgrids=\wannagridsx,
            ymajorgrids=\wannagridsy,
            grid style=dashed,
            xmode=log,
            log basis x=2,
            ymode=log,
            log basis y=10
        ]
    
        \addplot[
            color=blue,
            mark=square,
            ]
            coordinates {
            \coordinatesxya
            };
    
        \addplot[
            color=red,
            mark=pentagon*,
            ]
            coordinates {
            \coordinatesxyb
            };
        
        \addplot[
            color=green,
            mark=halfcircle*,
            ]
            coordinates {
            \coordinatesxyc
            };
        \legend{\legenda,\legendb,\legendc}
    
        \end{axis}
        \end{tikzpicture}
        \caption{Weak scaling, structured meshes. Comparison of CG iterations for the elliptic solvers (left)
and solution time for the elliptic and parabolic solvers (right) on CPU.}
    \end{figure}
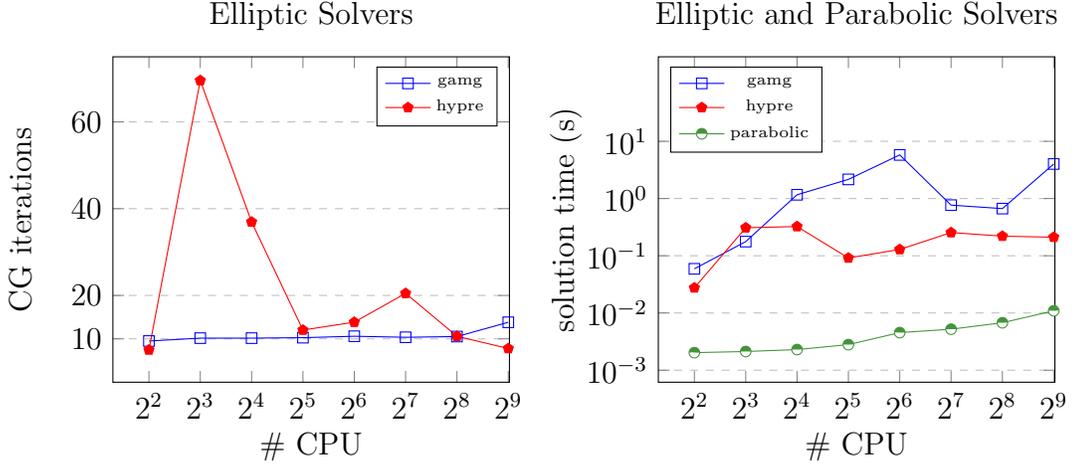

\subsection{Test 4 - simulation of the whole heartbeat}\label{heartbeat}
In this test, the simulation has been ran for a whole beat, namely 500 ms, on the structured mesh. \\
In particular, we have considered the maximum number of DOFs that can be handled by the memory (i.e. $8487168$, corresponding to a FEM mesh of $256\times 256\times 128$ elements) when Hypre is employed as preconditioner and solved the problem exploiting 4 GPUs, namely all the GPUs available for a single node. Then, we have done the same test for 32 CPUs, all the processors available on a single node for a $128\times 128\times 64$ mesh, yielding a total amount of 1073280 DOFs, namely the maximum number of DOFs that can be handled by the memory (240GB) when using GAMG as preconditioner. \\
We have also performed other tests on a $128\times 128\times 64$ mesh with 4 GPUs and 32 CPUs with Hypre in order to compare GAMG and Hypre performance on the same setup. In the following plots we have reported for each simulation timestep the number of iterations for the elliptic problem and the corresponding solution time, while in table \ref{beat} we have reported the same mean values benchmarked in the previous sections.

\begin{figure}[H]
    \includegraphics[width=\textwidth]{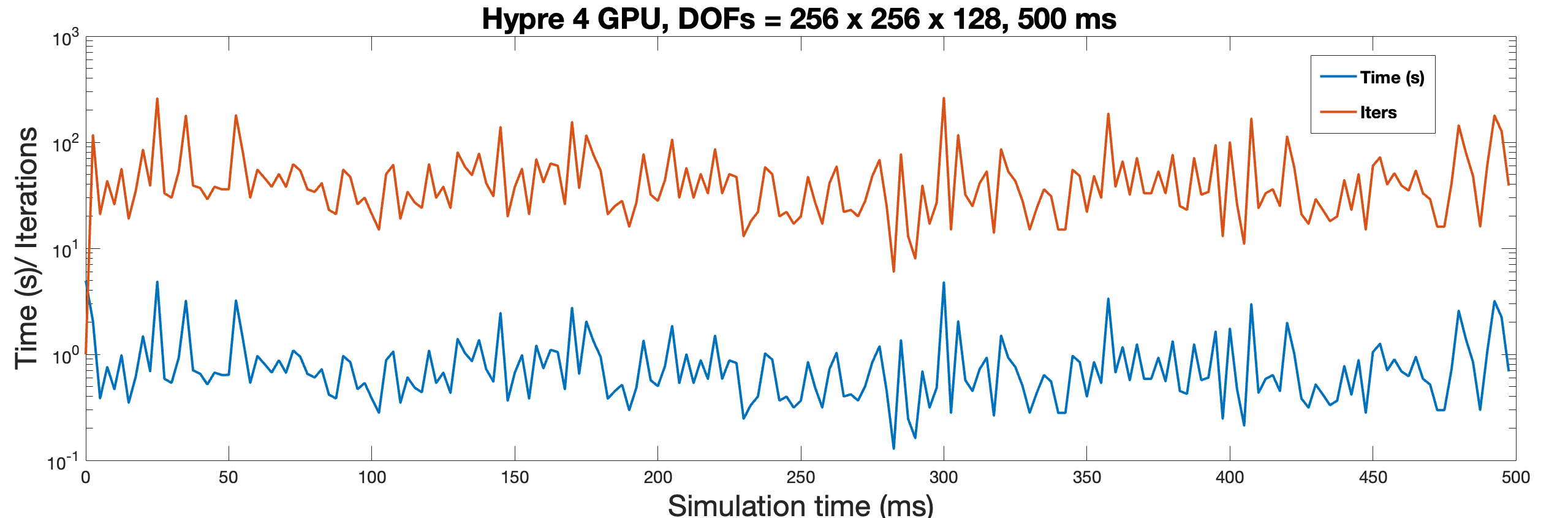}
    \caption{Simulation of the whole heartbeat on the structured meshes. Solution time and CG iterations for the elliptic system vs simulation time for Hypre on 4 GPUs with a mesh of $256\times 256\times 128$ elements.}
\end{figure}

\begin{figure}[H]
    \includegraphics[width=\textwidth]{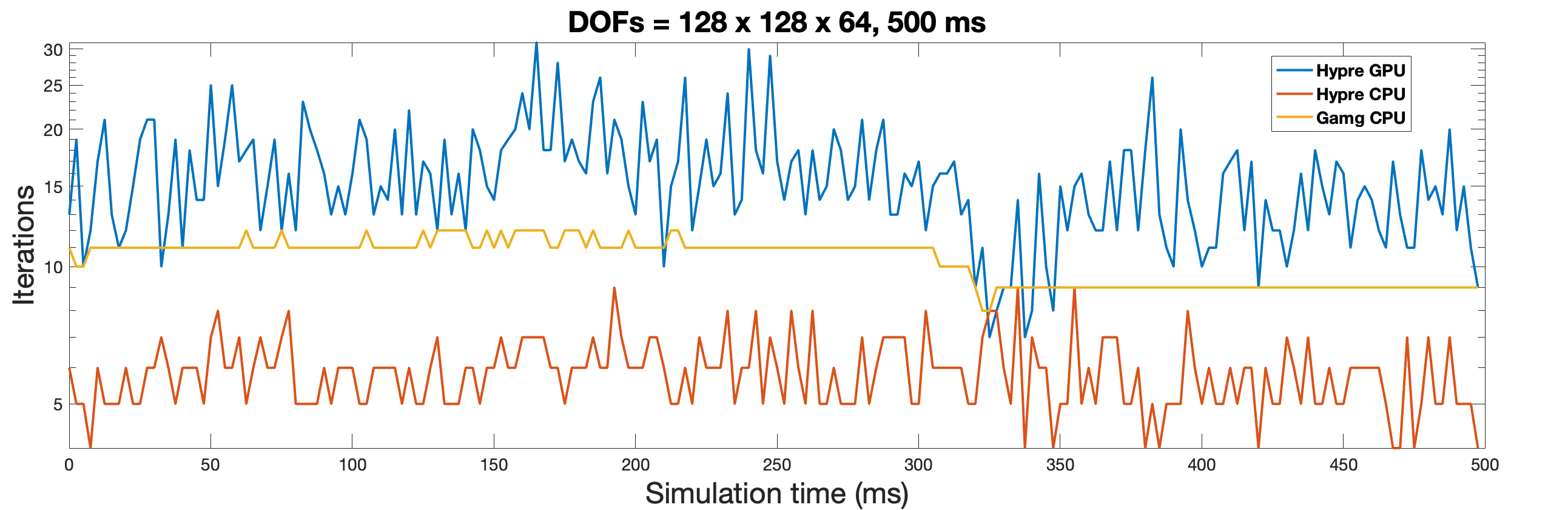}
    \caption{Simulation of the whole heartbeat on the structured meshes. Iterations for the elliptic system on a whole heartbeat simulation with a mesh of $128\times 128\times 64$ elements.}
\end{figure}

\begin{figure}[H]
    \includegraphics[width=\textwidth]{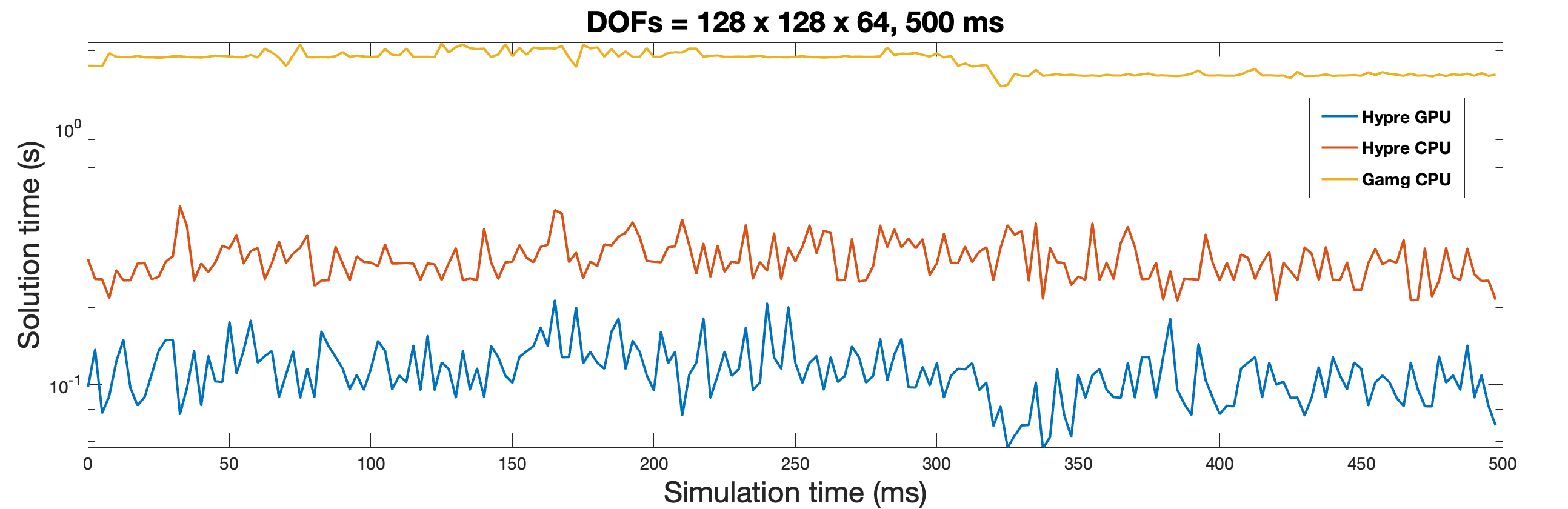}
    \caption{Simulation of the whole heartbeat on the structured meshes. Time to solve the elliptic system on a whole heartbeat simulation with a mesh of $128\times 128\times 64$ elements.}
\end{figure}

\begin{table}[H]
    \begin{tblr}{columns={c,font=\footnotesize}}
        \toprule
        Preconditioner & $\text{It}_{\text{parab, mean}}$ & $\text{It}_{\text{ellip, mean}}$ & $T_{\text{memb, mean}}$ (s) & $T_{\text{parab, mean}}$ (s) & $T_{\text{ellip, mean}}$ (s) \\
        \midrule
        Hypre 4 GPU & 26.00 & 43.62 & 4.7E-03 & 7.1E-02 & 0.80 \\
        \bottomrule
    \end{tblr}
    \caption{Mean times and iterations on a $256\times 256\times 128$ mesh for a whole heartbeat. Results obtained for Hypre on 4 GPUs.}
\end{table}

\begin{table}[H]
    \begin{tblr}{columns={c,font=\footnotesize}}
        \toprule
        Preconditioner & $\text{It}_{\text{parab, mean}}$ & $\text{It}_{\text{ellip, mean}}$ & $T_{\text{memb, mean}}$ (s) & $T_{\text{parab, mean}}$ (s) & $T_{\text{ellip, mean}}$ (s) \\
        \midrule
        GAMG 32 CPU & 6.53 & 10.18 & 1.7E-03 & 3.5E-02 & 1.81 \\
        Hypre 32 CPU & 6.53 & 4.44 & 4.7E-03 & 3.3E-02 & 0.24 \\
        Hypre 4 GPU & 26.98 & 15.22 & 7.7E-04 & 2.5E-02 & 0.11 \\
        \bottomrule
    \end{tblr}
    \caption{Mean times and iterations for Hypre and GAMG on a $128\times 128\times 64$ mesh for a whole heartbeat. On GPU the parabolic 
    problem is not preconditioned, thus we have an higher number of iterations.}
    \label{beat}
\end{table}

\section{Tests on Unstructured Mesh}\label{testunstr}
In this section we will comment the results of the parallel numerical tests
performed using the unstructured mesh as geometry for solving the Bidomain cardiac model. \\
First, in \textit{Test 1} we have studied the threshold parameters for the multigrid preconditioner
exploited for the elliptic equation. Then, in \textit{Test 2} we have performed a scaling test on 
three refinements of the geometry. Since the meshes considered are unstructured, a precise subdomain decomposition 
with an equal number of DOFs per worker has not been possible and the local DOFs have been calculated automatically
through the \texttt{PETSC\_DECIDE} option. 

\subsection{Test 1 - AMG Threshold Calibration}
For threshold tuning, we used the first mesh, U-mesh 1, for a total of 71450 DOFs, solving the membrane model with the GPU. 
The best results were obtained around 0.06-0.07 per GAMG, while for hypre, values between 0.5 and 0.6 were considered of particular interest.
The tests in this section were performed on a single node, in particular, if GPUs were used, all those available on the node (4) were used,
the same procedure was applied using CPUs: in this case, each node had a maximum of 32 physical cores. 

\begin{table}[H]
\centering
\begin{tblr}{columns={c,font=\footnotesize},vline{2}}
    \toprule
    Threshold & 0.25 & 0.3 & 0.4 & 0.5 & 0.6 & 0.7 \\
    \midrule
    $\text{It}_{\text{ellip, mean}}$ & 66.89 & 40.41 & 38.88 & 10.29 & 5.41 & 4.50 \\
    $T_{\text{ellip, mean}}$ (s) & 3.05e-01 & 1.77e-01 & 1.89e-01 & 5.24e-02 & 4.46e-02 & 4.46e-02 \\
    \bottomrule
\end{tblr}
\caption{AMG Threshold calibration, unstructured mesh (U-mesh 1, DOFs = 71450). Results for Hypre GPU solver. $\text{It}_{\text{ellip, mean}}$: average CG iterations per timestep for the elliptic solver. $T_{\text{ellip, mean}}$: average CG solution time (in s) per timestep for the elliptic solver.}
\label{tuhg}
\end{table}

\begin{table}[H]    
\centering
\begin{tblr}{columns={c,font=\footnotesize},vline{2}}
    \toprule
    Threshold & 0.25 & 0.3 & 0.4 & 0.5 & 0.6 & 0.7 \\
    \midrule
    $\text{It}_{\text{ellip, mean}}$ & 7.77 & 15.51 & 10.46 & 3.06 & 3.95 & 12.61 \\
    $T_{\text{ellip, mean}}$ (s) & 5.7E-02 & 1.0E-01 & 7.5E-02 & 3.0E-02 & 3.2E-02 & 7.6E-02 \\
    \bottomrule
\end{tblr}
\caption{AMG Threshold calibration, unstructured mesh (U-mesh 1, DOFs = 71450). Results for Hypre CPU solver. Same format as in Table \ref{tuhg}.}
\label{tuhc}
\end{table}

\begin{table}[H]
\centering
\begin{tblr}{columns={c,font=\footnotesize},vline{2}}
    \toprule
    Threshold & 0.0 & 0.01 & 0.02 & 0.03 & 0.04 & 0.05 & 0.06 & 0.07 \\
    \midrule
    $\text{It}_{\text{ellip, mean}}$ & 72.08 & 60.55 & 49.28 & 27.80 & 44.80 & 24.49 & 12.03 & 11.77 \\
    $T_{\text{ellip, mean}}$ (s) & 1.3E-01 & 1.1E-01 & 1.5E-01 & 5.8E-02 & 1.0E-01 & 6.6E-02 & 4.0E-02 & 4.3E-02 \\
    \bottomrule
\end{tblr}
\caption{AMG Threshold calibration, unstructured mesh (U-mesh 1, DOFs = 71450). Results for GAMG CPU solver. Same format as in Table \ref{tuhg}.}
\label{tugc}
\end{table}

\subsection{Test 2 - Strong Scaling}
In this test we have performed a strong scaling test, solving the problem on each of the three refinements considered.
We have performed tests on both CPU and GPU architectures, exploiting the algebraic multigrid implementations provided by
PETSc for preconditioning the elliptic system. Again, on CPU we have preconditioned the parabolic system using a block Jacobi preconditioner,
while on GPU the parabolic system is unpreconditioned. We notice in general good performance using GPUs instead of CPUs, with generally lower 
solution times expecially with an higher number of DOFs, which consent to fully exploit the potential of the GPU acceleration,
minimizing losses in time due mainly to synchronization overhead. We also noticed an out of memory error with U-mesh 3, the finest one
and more than 32 GPUs. This is probably a technical issue related to the handling of the copies for the acceleration device and it is still 
under investigation. However, up to 16 GPUs, performances are about one order of magnitude better than the CPU counterpart. \\
Figure \ref{fig:scaling_unst} reports the results of the strong scaling test performed on CPU. 

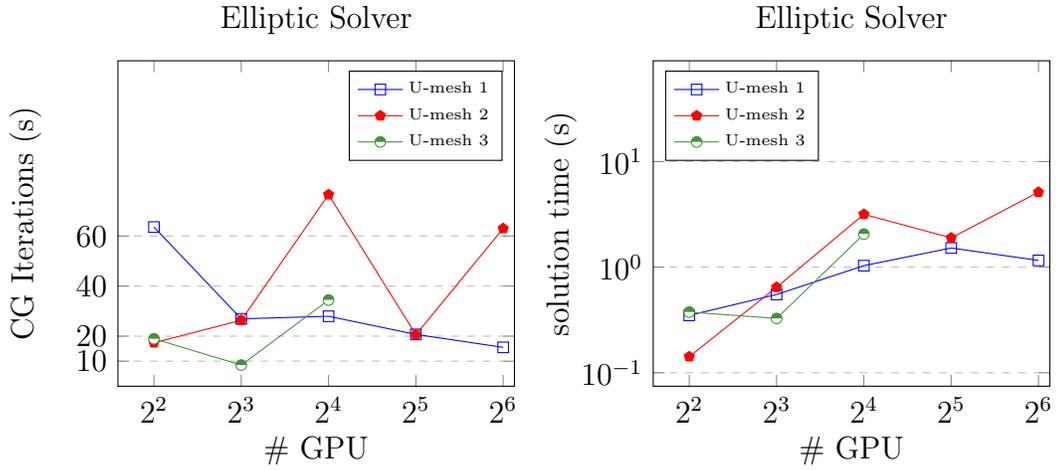
\begin{figure}[H] 
    \begin{tikzpicture}
        \def\titleax{Elliptic Solver}
        \def\titlex{\# GPU}
        \def\titley{CG Iterations (s)}
        \def\xmin{0}
        \def\xmax{70}
        \def\ymin{0}
        \def\ymax{130}
        \def\legpos{north east}
        \def\wannagridsx{false}
        \def\wannagridsy{true}
        
        \def\coordinatesxya{(4,6.3584e+01)(8,2.6911e+01)(16,2.7941e+01)(32,2.0693e+01)(64,1.5525e+01)}
        \def\legenda{U-mesh 1}
          
        \def\coordinatesxyb{(4,1.7347e+01)(8,2.6347e+01)(16,7.6525e+01)(32,2.0584e+01)(64,6.3010e+01)}
        \def\legendb{U-mesh 2}
    
        \def\coordinatesxyc{(4,1.8941e+01)(8,8.4455e+00)(16,3.4426e+01)}
        \def\legendc{U-mesh 3}
        
        \begin{axis}[
            width=0.5\textwidth,
            title={\titleax},
            xlabel={\titlex},
            ylabel={\titley},
            xmin=\xmin, xmax=\xmax,
            ymin=\ymin, ymax=\ymax,
            xtick={4,8,16,32,64},
            ytick={10,20,40,60},
            legend pos=\legpos,
            legend style={font=\tiny},
            xmajorgrids=\wannagridsx,
            ymajorgrids=\wannagridsy,
            grid style=dashed,
            xmode=log,
            log basis x=2
        ]
    
        \addplot[
            color=blue,
            mark=square,
            ]
            coordinates {
            \coordinatesxya
            };
    
        \addplot[
            color=red,
            mark=pentagon*,
            ]
            coordinates {
            \coordinatesxyb
            };
        \addplot[
            color=green,
            mark=halfcircle*,
            ]
            coordinates {
            \coordinatesxyc
            };
        
        \legend{\legenda,\legendb,\legendc}
    
        \end{axis}
    \end{tikzpicture}
    \begin{tikzpicture}
        \def\titleax{Elliptic Solver}
        \def\titlex{\# GPU}
        \def\titley{solution time (s)}
        \def\xmin{0}
        \def\xmax{70}
        \def\ymin{0}
        \def\ymax{90}
        \def\legpos{north west}
        \def\wannagridsx{false}
        \def\wannagridsy{true}
        \def\xticks{4,8,16,32,64}
        \def\yticks{1e-3,1e-2,1e-1,1,10}
        
        \def\coordinatesxya{(4,3.4916e-01)(8,5.5089e-01)(16,1.0331e+00)(32,1.5189e+00)(64,1.1581e+00)}
        \def\legenda{U-mesh 1}
        
        \def\coordinatesxyb{(4,1.4212e-01)(8,6.4500e-01)(16,3.1671e+00)(32,1.8942e+00)(64,5.1272e+00)}
        \def\legendb{U-mesh 2}
               
        \def\coordinatesxyc{(4,3.7547e-01)(8,3.2656e-01)(16,2.0549e+00)}
        \def\legendc{U-mesh 3}
    
        \begin{axis}[
            width=0.5\textwidth,
            title={\titleax},
            xlabel={\titlex},
            ylabel={\titley},
            xmin=\xmin, xmax=\xmax,
            ymin=\ymin, ymax=\ymax,
            xtick={4,8,16,32,64,128,256,512},
            ytick={1e-3,1e-2,1e-1,1,10},
            legend pos=\legpos,
            legend style={font=\tiny},
            xmajorgrids=\wannagridsx,
            ymajorgrids=\wannagridsy,
            grid style=dashed,
            xmode=log,
            log basis x=2,
            ymode=log,
            log basis y=10
        ]
    
        \addplot[
            color=blue,
            mark=square,
            ]
            coordinates {
            \coordinatesxya
            };
    
        \addplot[
            color=red,
            mark=pentagon*,
            ]
            coordinates {
            \coordinatesxyb
            };
        
        \addplot[
            color=green,
            mark=halfcircle*,
            ]
            coordinates {
            \coordinatesxyc
            };
        \legend{\legenda,\legendb,\legendc}
    
        \end{axis}
        \end{tikzpicture}
        \caption{Strong scaling, unstructured meshes. Comparison of CG iterations (left) and solution time (right) for the elliptic solver on GPU (Hypre BoomerAMG). U-mesh 3 goes out of memory with 32 and 64 GPUs.}
    \end{figure}

\begin{figure}[H] 
    \begin{tikzpicture}
        \def\titleax{Elliptic Solver}
        \def\titlex{\# CPU}
        \def\titley{CG iterations}
        \def\xmin{0}
        \def\xmax{520}
        \def\ymin{0}
        \def\ymax{70}
        \def\legpos{north east}
        \def\wannagridsx{false}
        \def\wannagridsy{true}
        
        \def\coordinatesxya{(4,1.1802e+01)(8,8.4059e+00)(16,8.5248e+00)(32,8.8614e+00)(64,1.6079e+01)(128,9.0693e+00)(256,1.1426e+01)(512,6.8317e+00)}
        \def\legenda{U-mesh 1}
          
        \def\coordinatesxyb{(4,8.3861e+00)(8,3.4772e+01)(16,4.8584e+01)(32,1.7624e+01)(64,2.8416e+01)(128,8.6337e+00)(256,2.2762e+01)(512,7.5644e+00)}
        \def\legendb{U-mesh 2}
    
        \def\coordinatesxyc{(4,3.3653e+01)(8,3.3891e+01)(16,2.3990e+01)(32,2.3129e+01)(64,3.6267e+01)(128,3.6743e+01)(256,3.5515e+01)(512,6.1485e+00)}
        \def\legendc{U-mesh 3}
        
        \begin{axis}[
            width=0.5\textwidth,
            title={\titleax},
            xlabel={\titlex},
            ylabel={\titley},
            xmin=\xmin, xmax=\xmax,
            ymin=\ymin, ymax=\ymax,
            xtick={4,8,16,32,64,128,256,512},
            ytick={10,30,50,60},
            legend pos=\legpos,
            legend style={font=\tiny},
            xmajorgrids=\wannagridsx,
            ymajorgrids=\wannagridsy,
            grid style=dashed,
            xmode=log,
            log basis x=2
        ]
    
        \addplot[
            color=blue,
            mark=square,
            ]
            coordinates {
            \coordinatesxya
            };
    
        \addplot[
            color=red,
            mark=pentagon*,
            ]
            coordinates {
            \coordinatesxyb
            };
        \addplot[
            color=green,
            mark=halfcircle*,
            ]
            coordinates {
            \coordinatesxyc
            };
        
        \legend{\legenda,\legendb,\legendc}
    
        \end{axis}
    \end{tikzpicture}
    \begin{tikzpicture}
        \def\titleax{Elliptic Solver}
        \def\titlex{\# CPU}
        \def\titley{solution time (s)}
        \def\xmin{0}
        \def\xmax{520}
        \def\ymin{0}
        \def\ymax{110}
        \def\legpos{north east}
        \def\wannagridsx{false}
        \def\wannagridsy{true}
        \def\xticks{4,8,16,32,64,128,256,512}
        \def\yticks{1e-2,1e-1,1,10,100}
        
        \def\coordinatesxya{(4,2.0653e-01)(8,9.3492e-02)(16,6.3081e-02)(32,6.5680e-02)(64,1.4327e-01)(128,1.7095e-01)(256,5.3590e-01)(512,4.6408e-01)}
        \def\legenda{U-mesh 1}
        
        \def\coordinatesxyb{(4,1.2902e+00)(8,2.6865e+00)(16,2.0666e+00)(32,5.3947e-01)(64,6.1860e-01)(128,2.8729e-01)(256,1.4559e+00)(512,9.1906e-01)}
        \def\legendb{U-mesh 2}
               
        \def\coordinatesxyc{(4,4.0513e+01)(8,2.1719e+01)(16,8.1705e+00)(32,5.3632e+00)(64,4.3110e+00)(128,3.0026e+00)(256,4.1072e+00)(512,1.2034e+00)}
        \def\legendc{U-mesh 3}
    
        \begin{axis}[
            width=0.5\textwidth,
            title={\titleax},
            xlabel={\titlex},
            ylabel={\titley},
            xmin=\xmin, xmax=\xmax,
            ymin=\ymin, ymax=\ymax,
            xtick={4,8,16,32,64,128,256,512},
            ytick={1e-3,1e-2,1e-1,1,10},
            legend pos=\legpos,
            legend style={font=\tiny},
            xmajorgrids=\wannagridsx,
            ymajorgrids=\wannagridsy,
            grid style=dashed,
            xmode=log,
            log basis x=2,
            ymode=log,
            log basis y=10
        ]
    
        \addplot[
            color=blue,
            mark=square,
            ]
            coordinates {
            \coordinatesxya
            };
    
        \addplot[
            color=red,
            mark=pentagon*,
            ]
            coordinates {
            \coordinatesxyb
            };
        
        \addplot[
            color=green,
            mark=halfcircle*,
            ]
            coordinates {
            \coordinatesxyc
            };
        \legend{\legenda,\legendb,\legendc}
    
        \end{axis}
        \end{tikzpicture}
        \caption{
        Strong scaling, unstructured meshes. Comparison of CG iterations (left) and solution time (right) for the elliptic solver on CPU (Hypre BoomerAMG).
        }\label{fig:scaling_unst}
    \end{figure}

\begin{figure}[H] 
    \begin{tikzpicture}
        \def\titleax{Elliptic Solver}
        \def\titlex{\# CPU}
        \def\titley{Iterations}
        \def\xmin{0}
        \def\xmax{520}
        \def\ymin{5}
        \def\ymax{30}
        \def\legpos{north east}
        \def\wannagridsx{false}
        \def\wannagridsy{true}
        
        \def\coordinatesxya{(4,1.1238e+01)(8,1.1040e+01)(16,1.1287e+01)(32,1.1772e+01)(64,1.1624e+01)(128,1.1010e+01)(256,1.1079e+01)(512,1.0970e+01)}
        \def\legenda{U-mesh 1}
          
        \def\coordinatesxyb{(4,1.3871e+01)(8,1.3990e+01)(16,1.3634e+01)(32,1.3941e+01)(64,1.3842e+01)(128,1.4089e+01)(256,1.4139e+01)(512,1.3861e+01)}
        \def\legendb{U-mesh 2}
    
        \def\coordinatesxyc{(4,1.7149e+01)(8,1.7267e+01)(16,1.7772e+01)(32,1.7713e+01)(64,1.7812e+01)(128,1.7990e+01)(256,1.7653e+01)(512,1.7535e+01)}
        \def\legendc{U-mesh 3}
        
        \begin{axis}[
            width=0.5\textwidth,
            title={\titleax},
            xlabel={\titlex},
            ylabel={\titley},
            xmin=\xmin, xmax=\xmax,
            ymin=\ymin, ymax=\ymax,
            xtick={4,8,16,32,64,128,256,512},
            ytick={10,15,20},
            legend pos=\legpos,
            legend style={font=\tiny},
            xmajorgrids=\wannagridsx,
            ymajorgrids=\wannagridsy,
            grid style=dashed,
            xmode=log,
            log basis x=2
        ]
    
        \addplot[
            color=blue,
            mark=square,
            ]
            coordinates {
            \coordinatesxya
            };
    
        \addplot[
            color=red,
            mark=pentagon*,
            ]
            coordinates {
            \coordinatesxyb
            };
        \addplot[
            color=green,
            mark=halfcircle*,
            ]
            coordinates {
            \coordinatesxyc
            };
        
        \legend{\legenda,\legendb,\legendc}
    
        \end{axis}
    \end{tikzpicture}
    \begin{tikzpicture}
        \def\titleax{Elliptic Solver}
        \def\titlex{\# CPU}
        \def\titley{time (s)}
        \def\xmin{0}
        \def\xmax{520}
        \def\ymin{0}
        \def\ymax{110}
        \def\legpos{north east}
        \def\wannagridsx{false}
        \def\wannagridsy{true}
        \def\xticks{4,8,16,32,64,128,256,512}
        \def\yticks{1e-2,1e-1,1,10,100}
        
        \def\coordinatesxya{(4,1.2103e-01)(8,7.4674e-02)(16,4.7473e-02)(32,3.9968e-02)(64,4.4045e-02)(128,7.4162e-02)(256,1.0698e-01)(512,1.4178e-01)}
        \def\legenda{U-mesh 1}
        
        \def\coordinatesxyb{(4,1.0247e+00)(8,5.7752e-01)(16,3.2948e-01)(32,2.3511e-01)(64,1.7391e-01)(128,1.9148e-01)(256,2.9924e-01)(512,4.6747e-01)}
        \def\legendb{U-mesh 2}
               
        \def\coordinatesxyc{(4,1.3494e+01)(8,7.4201e+00)(16,3.9965e+00)(32,2.4760e+00)(64,1.5094e+00)(128,1.0505e+00)(256,1.1795e+00)(512,1.3474e+00)}
        \def\legendc{U-mesh 3}
    
        \begin{axis}[
            width=0.5\textwidth,
            title={\titleax},
            xlabel={\titlex},
            ylabel={\titley},
            xmin=\xmin, xmax=\xmax,
            ymin=\ymin, ymax=\ymax,
            xtick={4,8,16,32,64,128,256,512},
            ytick={1e-3,1e-2,1e-1,1,10},
            legend pos=\legpos,
            legend style={font=\tiny},
            xmajorgrids=\wannagridsx,
            ymajorgrids=\wannagridsy,
            grid style=dashed,
            xmode=log,
            log basis x=2,
            ymode=log,
            log basis y=10
        ]
    
        \addplot[
            color=blue,
            mark=square,
            ]
            coordinates {
            \coordinatesxya
            };
    
        \addplot[
            color=red,
            mark=pentagon*,
            ]
            coordinates {
            \coordinatesxyb
            };
        
        \addplot[
            color=green,
            mark=halfcircle*,
            ]
            coordinates {
            \coordinatesxyc
            };
        \legend{\legenda,\legendb,\legendc}
    
        \end{axis}
        \end{tikzpicture}
        \caption{Strong scaling, unstructured meshes. Comparison of CG iterations (left) and solution time (right) for the elliptic solver on CPU (GAMG).}
    \end{figure}
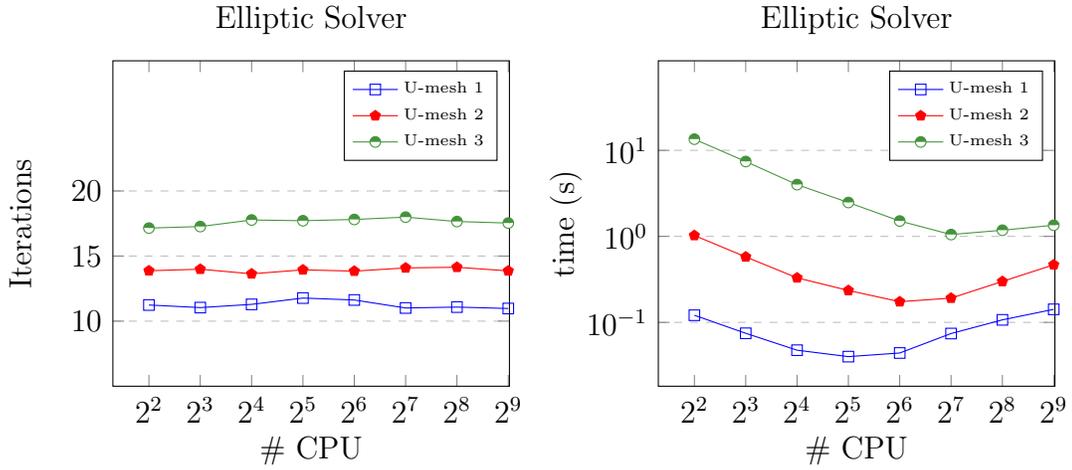
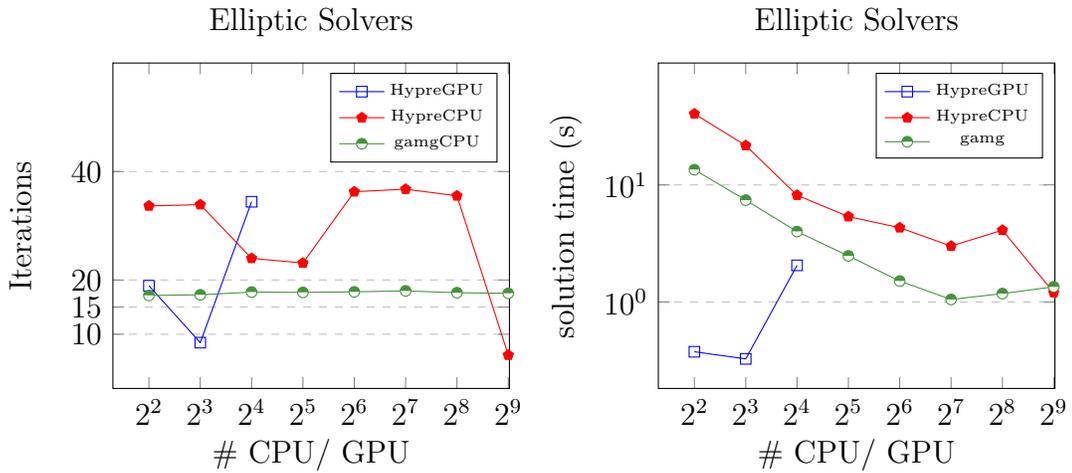
\begin{figure}[H] 
    \begin{tikzpicture}
        \def\titleax{Elliptic Solvers}
        \def\titlex{\# CPU/ GPU}
        \def\titley{Iterations}
        \def\xmin{0}
        \def\xmax{520}
        \def\ymin{0}
        \def\ymax{60}
        \def\legpos{north east}
        \def\wannagridsx{false}
        \def\wannagridsy{true}

        \def\coordinatesxya{(4,1.8941e+01)(8,8.4455e+00)(16,3.4426e+01)}
        \def\legenda{HypreGPU}
        
        \def\coordinatesxyb{(4,3.3653e+01)(8,3.3891e+01)(16,2.3990e+01)(32,2.3129e+01)(64,3.6267e+01)(128,3.6743e+01)(256,3.5515e+01)(512,6.1485e+00)}
        \def\legendb{HypreCPU}
        
        \def\coordinatesxyc{(4,1.7149e+01)(8,1.7267e+01)(16,1.7772e+01)(32,1.7713e+01)(64,1.7812e+01)(128,1.7990e+01)(256,1.7653e+01)(512,1.7535e+01)}
        \def\legendc{gamgCPU}
        
        \begin{axis}[
            width=0.5\textwidth,
            title={\titleax},
            xlabel={\titlex},
            ylabel={\titley},
            xmin=\xmin, xmax=\xmax,
            ymin=\ymin, ymax=\ymax,
            xtick={4,8,16,32,64,128,256,512},
            ytick={10,15,20,40},
            legend pos=\legpos,
            legend style={font=\tiny},
            xmajorgrids=\wannagridsx,
            ymajorgrids=\wannagridsy,
            grid style=dashed,
            xmode=log,
            log basis x=2
        ]
    
        \addplot[
            color=blue,
            mark=square,
            ]
            coordinates {
            \coordinatesxya
            };
    
        \addplot[
            color=red,
            mark=pentagon*,
            ]
            coordinates {
            \coordinatesxyb
            };

        \addplot[
            color=green,
            mark=halfcircle*,
            ]
            coordinates {
            \coordinatesxyc
            };
        
        \legend{\legenda,\legendb, \legendc}
    
        \end{axis}
    \end{tikzpicture}
    \begin{tikzpicture}
        \def\titleax{Elliptic Solvers}
        \def\titlex{\# CPU/ GPU}
        \def\titley{solution time (s)}
        \def\xmin{0}
        \def\xmax{520}
        \def\ymin{0}
        \def\ymax{110}
        \def\legpos{north east}
        \def\wannagridsx{false}
        \def\wannagridsy{true}
        \def\xticks{4,8,16,32,64,128,256,512}
        \def\yticks{1e-2,1e-1,1,10,100}

        \def\coordinatesxya{(4,3.7547e-01)(8,3.2656e-01)(16,2.0549e+00)}
        \def\legenda{HypreGPU}
        
        \def\coordinatesxyb{(4,4.0513e+01)(8,2.1719e+01)(16,8.1705e+00)(32,5.3632e+00)(64,4.3110e+00)(128,3.0026e+00)(256,4.1072e+00)(512,1.2034e+00)}
        \def\legendb{HypreCPU}
        
        \def\coordinatesxyc{(4,1.3494e+01)(8,7.4201e+00)(16,3.9965e+00)(32,2.4760e+00)(64,1.5094e+00)(128,1.0505e+00)(256,1.1795e+00)(512,1.3474e+00)}
        \def\legendc{gamg}
    
        \begin{axis}[
            width=0.5\textwidth,
            title={\titleax},
            xlabel={\titlex},
            ylabel={\titley},
            xmin=\xmin, xmax=\xmax,
            ymin=\ymin, ymax=\ymax,
            xtick={4,8,16,32,64,128,256,512},
            ytick={1e-3,1e-2,1e-1,1,10},
            legend pos=\legpos,
            legend style={font=\tiny},
            xmajorgrids=\wannagridsx,
            ymajorgrids=\wannagridsy,
            grid style=dashed,
            xmode=log,
            log basis x=2,
            ymode=log,
            log basis y=10
        ]
    
        \addplot[
            color=blue,
            mark=square,
            ]
            coordinates {
            \coordinatesxya
            };
    
        \addplot[
            color=red,
            mark=pentagon*,
            ]
            coordinates {
            \coordinatesxyb
            };

        \addplot[
            color=green,
            mark=halfcircle*,
            ]
            coordinates {
            \coordinatesxyc
            };
            
        \legend{\legenda,\legendb,\legendc}
    
        \end{axis}
        \end{tikzpicture}
        \caption{Strong scaling, unstructured meshes (U-mesh3). Comparison of CG iterations (left) and solution time (right) for the elliptic solver on CPU (GAMG and Hypre BoomerAMG) and on GPU (Hypre BoomerAMG).}
    \end{figure}

\begin{figure}[H]
    \centering
    \begin{tblr}{colspec = {X[c] X[c]}, colsep = 2pt}
    \includegraphics[width=0.5\linewidth]{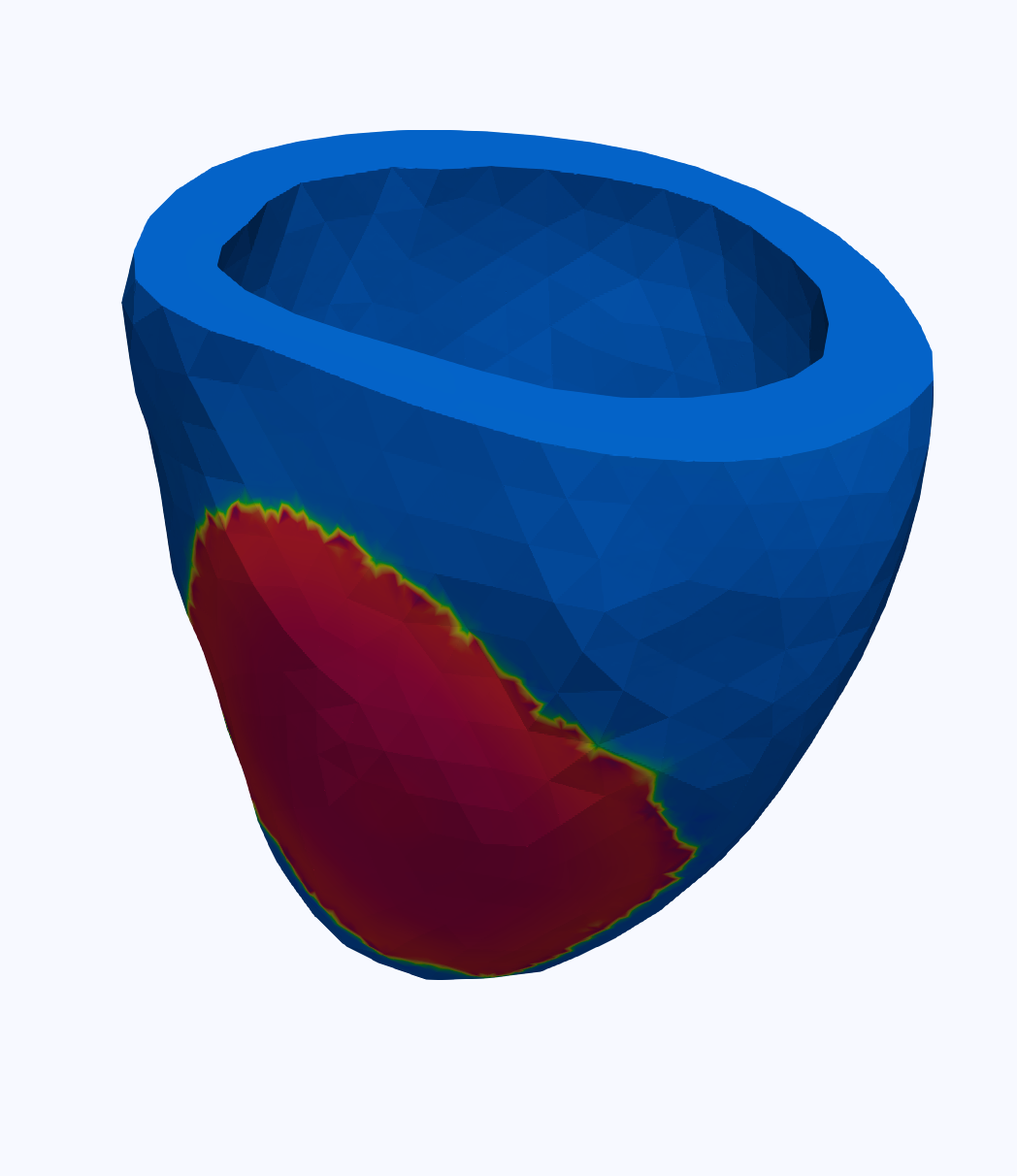} & \includegraphics[width=0.5\linewidth]{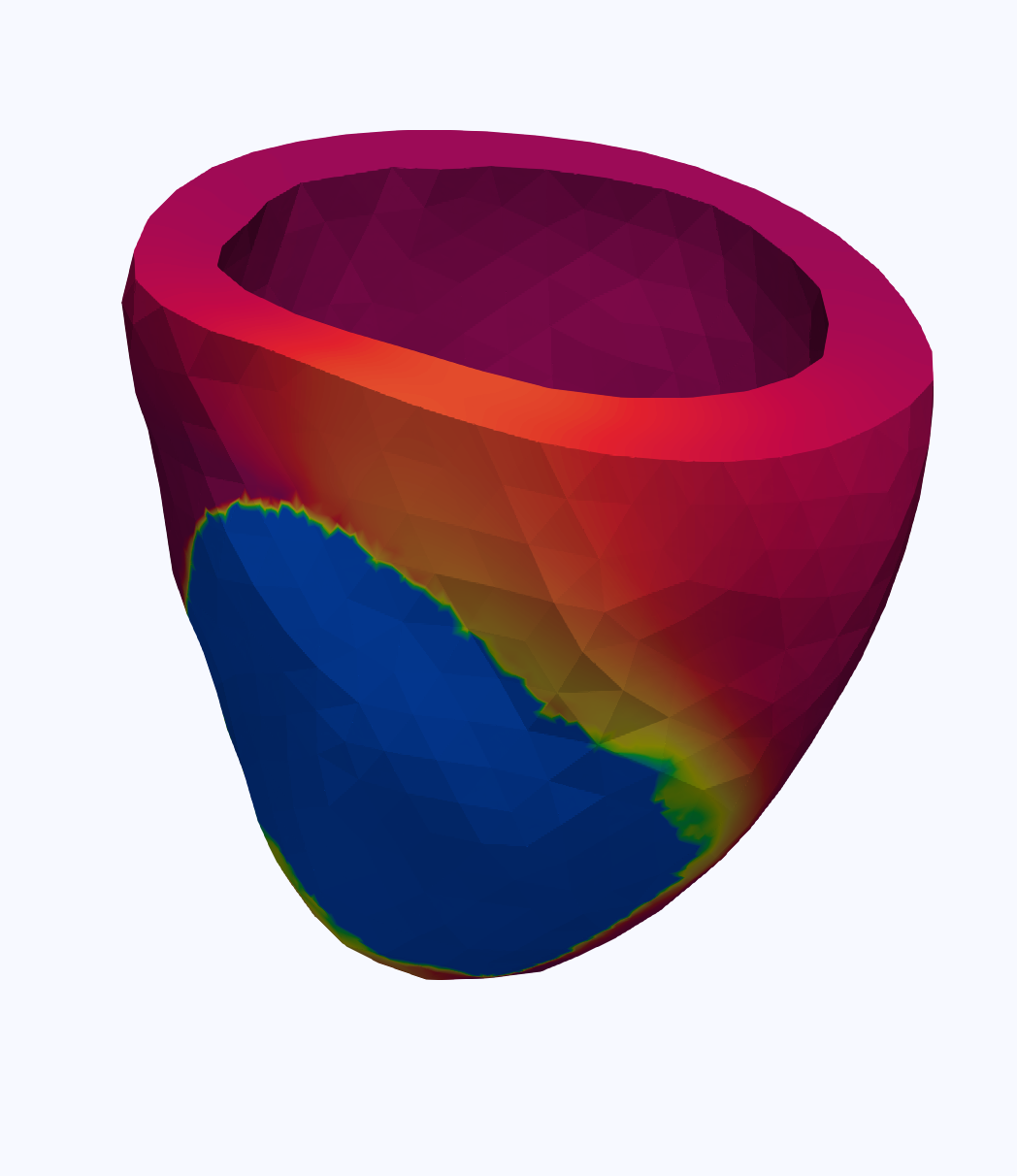} \\
    \includegraphics[width=0.5\linewidth]{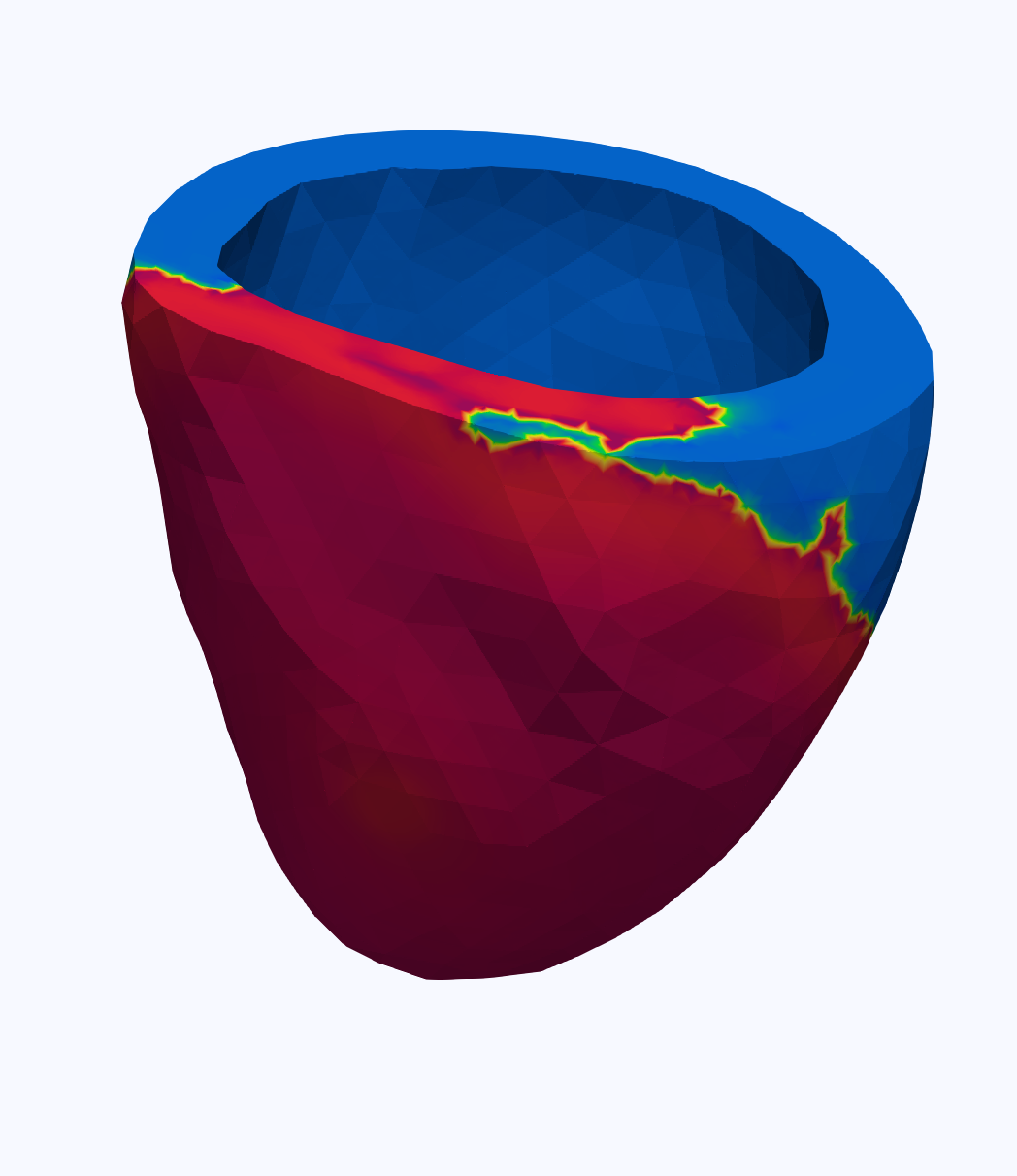} & \includegraphics[width=0.5\linewidth]{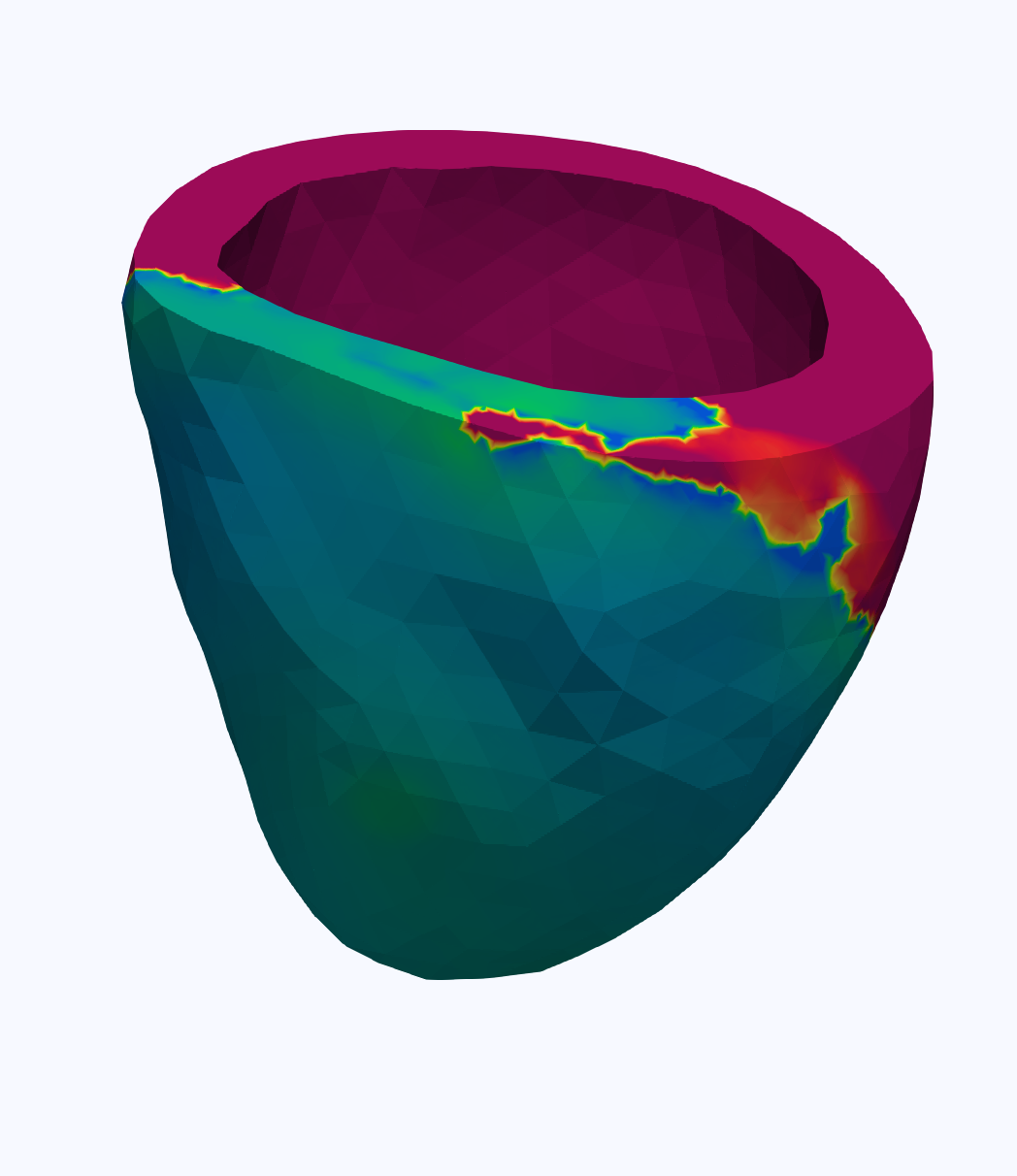} \\
    \includegraphics[width=0.5\linewidth]{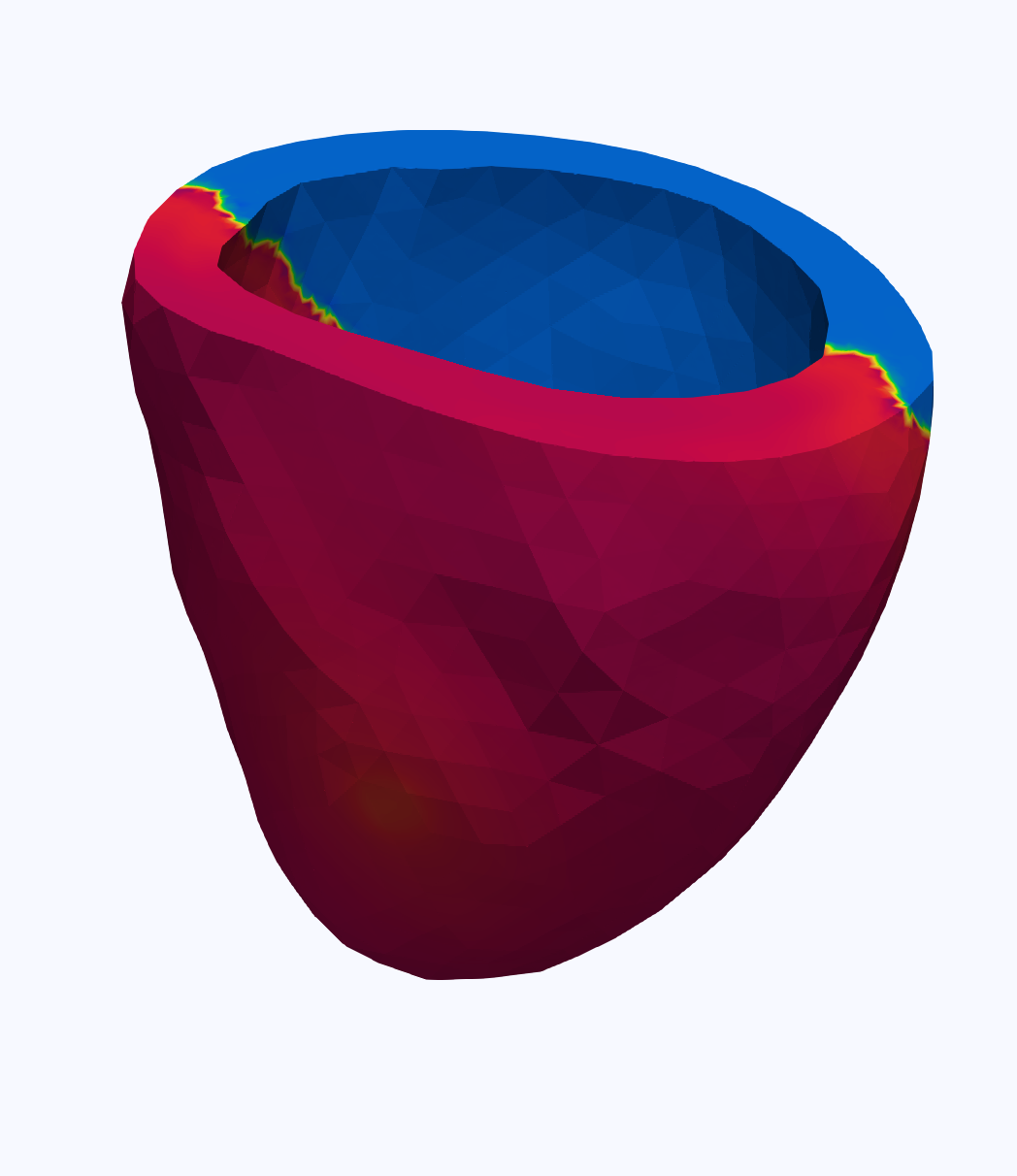} & \includegraphics[width=0.5\linewidth]{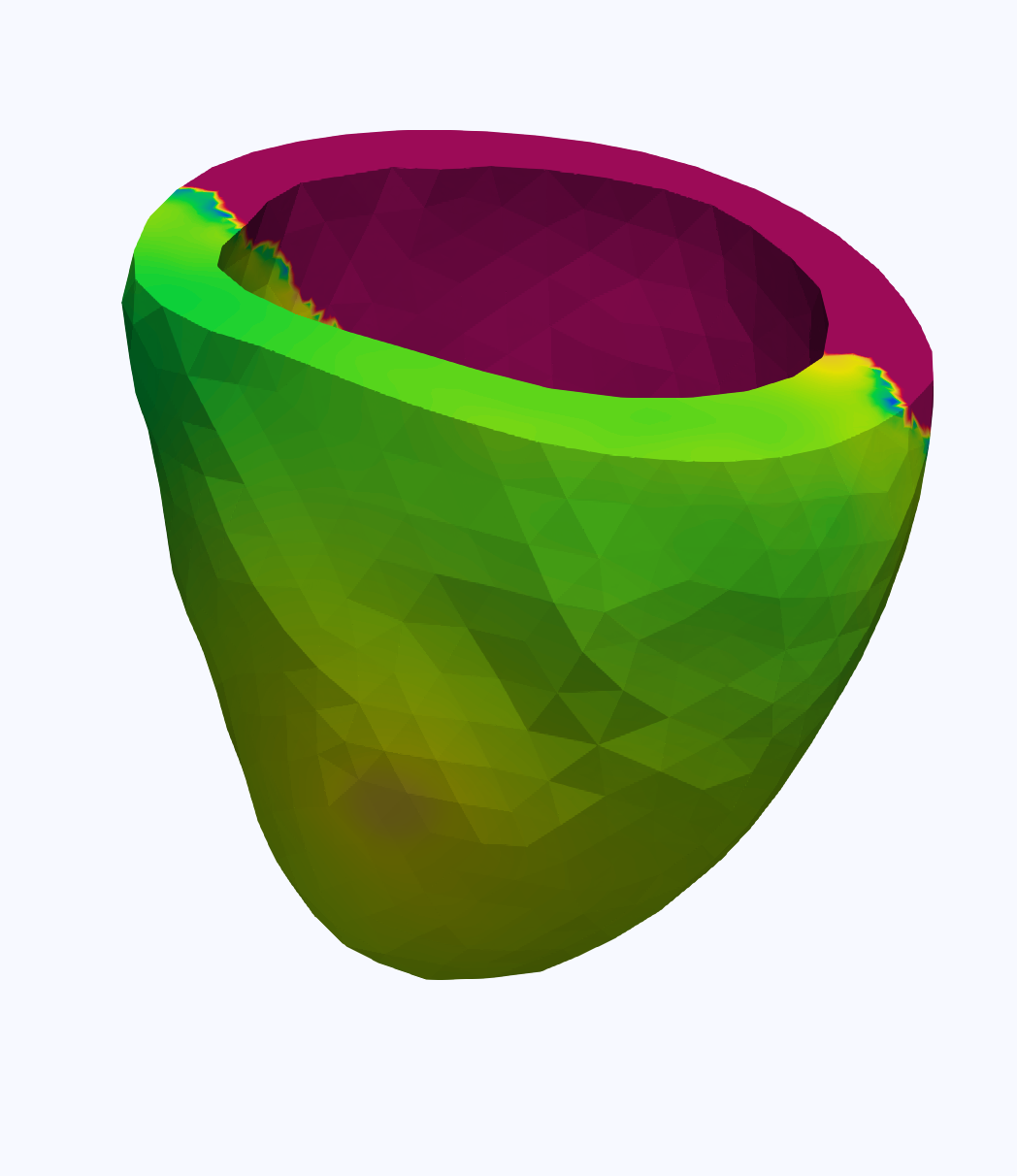} \\
    \includegraphics[width=0.5\linewidth]{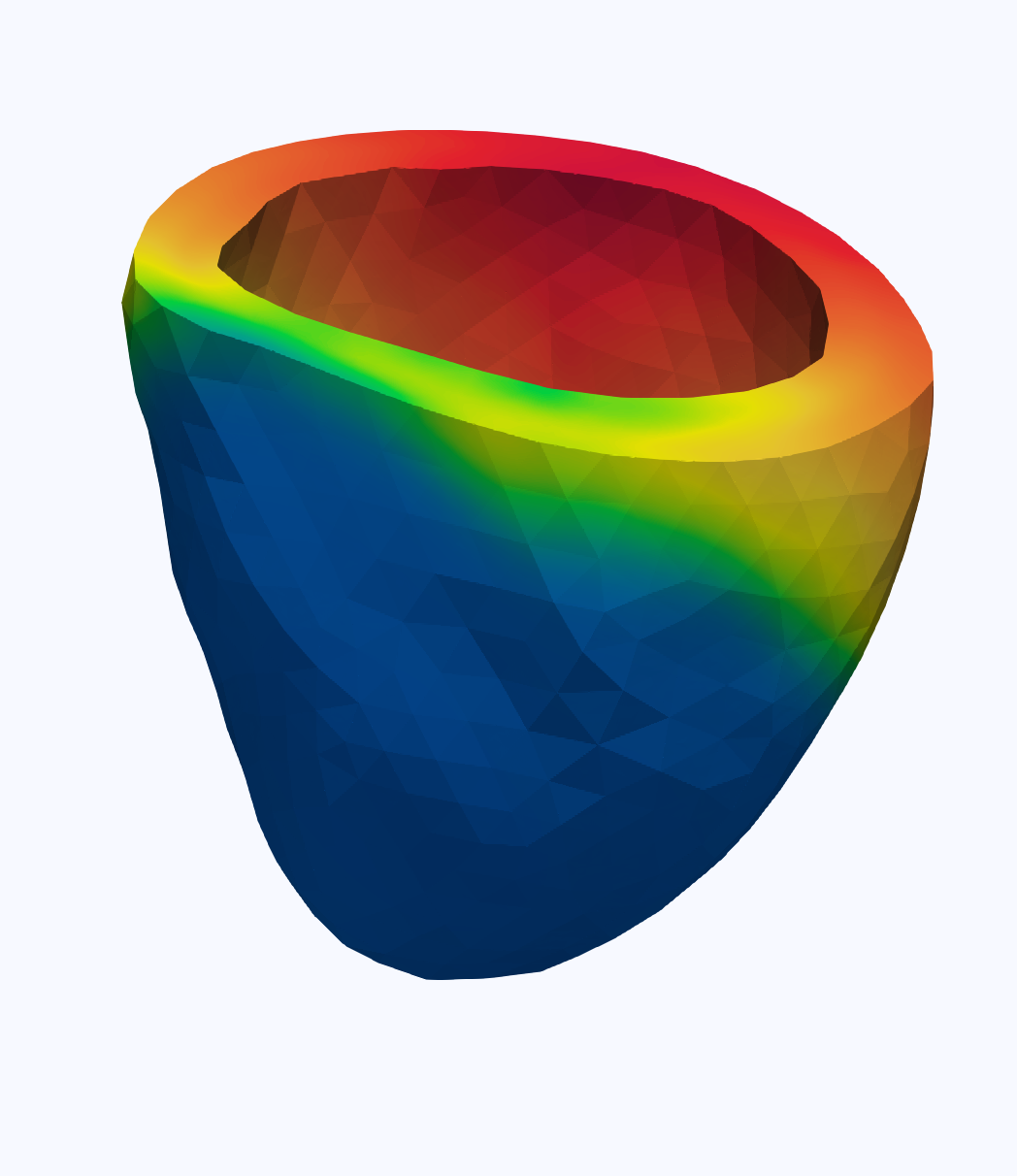} & \includegraphics[width=0.5\linewidth]{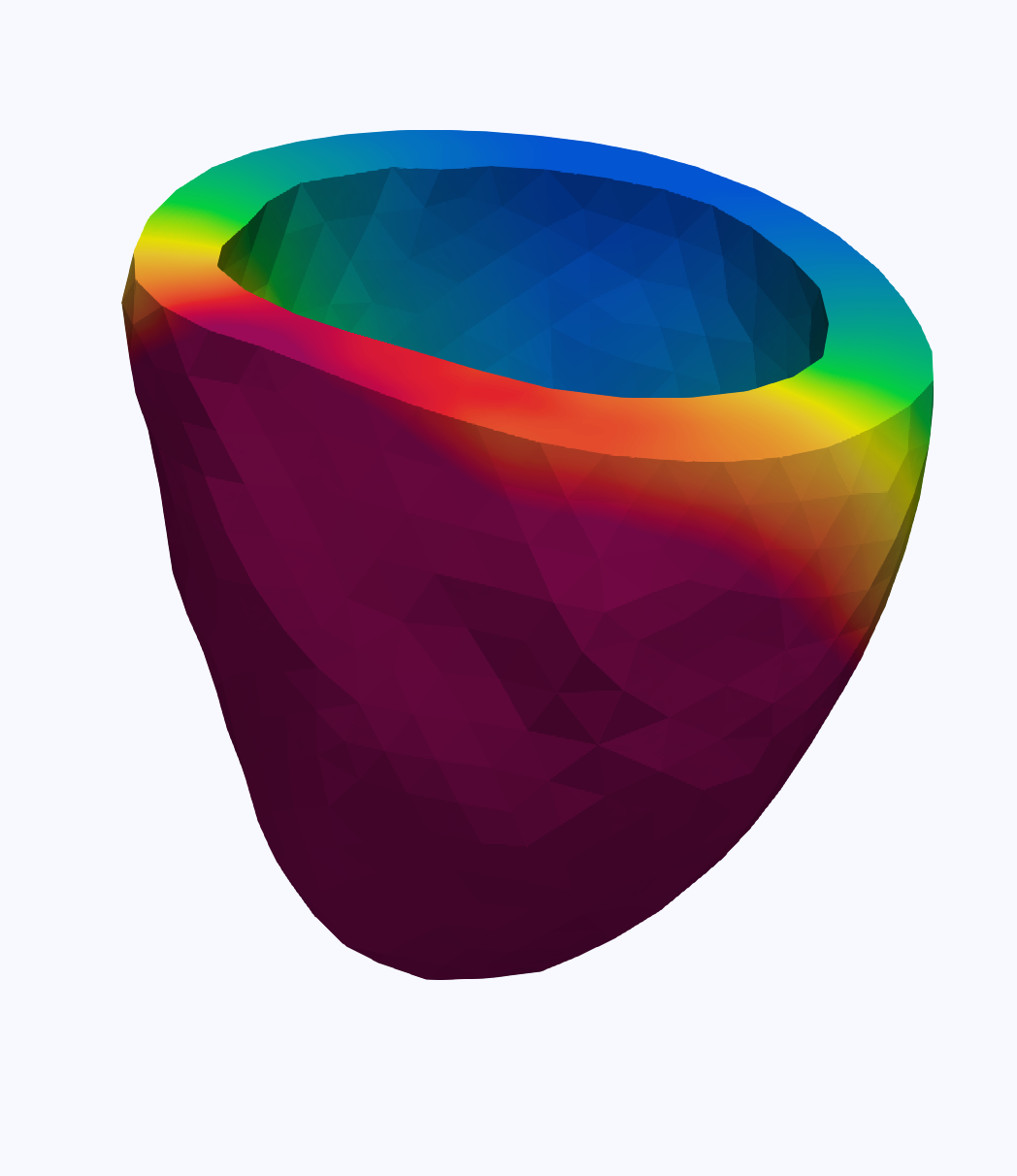} \\
    \includegraphics[width=0.5\linewidth]{color_scale_uu_crump.png} & \includegraphics[width=0.5\linewidth]{color_scale_ue_crump.png} \\
    \end{tblr}
\caption{Transmembrane potential (left) and extracellular potential (right) snapshots on the epicardial surface, represented by the U-mesh. The values of the displayed map are in mV}
\end{figure}

\section{Conclusions}
In this work we have performed scalability tests for the Bidomain cardiac model, focusing in particular on the performance
of the algebraic multigrid implementations provided by PETSc (GAMG and Hypre) for preconditioning the linear system
arising from the discretization of the elliptic system in the parabolic-elliptic formulation of the model. \\
Tests were performed on CPU and GPU on two different kind of meshes: a structured one and a more realistic unstructured one, representing 
a cardiac ventricle. Overall results have shown a general better scaling properties (expecially strong scaling) when the problem is solved on CPU, even if
the absolute best performance was obtained when solving the problem on GPU. In particular, in case of structured meshes we observed that the solution of the elliptic problem on GPU is 3.6 times faster than the CPU counterpart, whereas in case of unstructured meshes the solution of the elliptic problem is 3.2 times faster than on CPU.
From the tests presented we have also confirmed the scalability properties of AMG used as preconditioner for a more complex problem than the one considered in \cite{cle00} and we have provided an extension of the results in \cite{nei12} with more GPUs, different numerical schemes and an updated version of the software employed.  

\appendix
\section{PETSc Options}\label{opts}
In this section we report the PETSc options set for our numerical parallel tests. \\
For Hypre BoomerAMG, we have considered default options except for the following:

\begin{lstlisting}[basicstyle=\footnotesize]
    -dm_mat_type $MAT_TYPE                    # mpiaij or aijcusparse
    -dm_vec_type $VEC_TYPE                    # mpi or cuda
    -pc_hypre_boomeramg_strong_threshold $THR # Threshold value
\end{lstlisting}

For GAMG instead, we have considered the following \textit{non-default} options
\begin{lstlisting}[basicstyle=\footnotesize]
    -dm_mat_type $MAT_TYPE                            # mpiaij or aijcusparse
    -dm_vec_type $VEC_TYPE                            # mpi or cuda
    -pc_GAMG_type agg
    -pc_GAMG_agg_nsmooths 1
    -pc_GAMG_coarse_eq_limit 100
    -pc_GAMG_reuse_interpolation
    -pc_GAMG_square_graph 1
    -pc_GAMG_threshold $THR                           # Threshold value
    -mg_levels_ksp_max_it 2
    -mg_levels_ksp_type chebyshev
    -mg_levels_esteig_ksp_type cg
    -mg_levels_esteig_ksp_max_it 10
    -mg_levels_ksp_chebyshev_esteig 0,0.05,0,1.05
    -mg_levels_pc_type jacobi
\end{lstlisting}

\section*{Acknowledgements}

We would like to thank Marco Fedele from Politecnico di Milano for providing us the patient-specific LV mesh and Luca F. Pavarino
from the University of Pavia for many helpful discussions and suggestions.

This work was part of the MICROCARD project. This project has received funding from the European High-Performance Computing Joint Undertaking EuroHPC (JU) under grant agreement No 955495. The JU receives support from the European Union’s Horizon 2020 research and innovation programme and France, Italy, Germany, Austria, Norway, Switzerland.
E. Centofanti and S. Scacchi have been supported by grants of INdAM--GNCS. S. Scacchi has been supported by MIUR grants PRIN 2017AXL54F$\_$003 and PRIN 202232A8AN$\_$003.

\end{document}